\documentclass[reqno]{amsart}
\usepackage{graphicx}
\usepackage{latexsym}
\usepackage{amssymb}
\usepackage{amsmath}
\usepackage{layout}

\newtheorem{theorem}{Theorem}[section]

\newtheorem{lemma}[theorem]{Lemma}

\newtheorem{example}[theorem]{Example}
\newtheorem{remark}[theorem]{Remark}
\newtheorem{definition}[theorem]{Definition}

\def\real{{\mathord{{\rm I\kern-2.8pt R}}}}        
\def\inte{{\mathord{{\rm I\kern-2.8pt N}}}}

\def\sZZ{{\rm Z\kern-2.8ptem{}Z}}

\def\z{{\mathchoice
  {\sZZ}
  {\sZZ}
  {\rm Z\kern-0.30em{}Z}
  {\rm Z\kern-0.25em{}Z} }}
\def\sQQ{{\kern 0.27em \vrule height1.45ex width0.03em depth0em
          \kern-0.30em \rm Q}}
\def\qu{{\mathchoice
    {\sQQ}
    {\sQQ}
  {\kern 0.225em \vrule height1.05ex width0.025em depth0em \kern-0.25em \rm Q}
  {\kern 0.180em \vrule height0.78ex width0.020em depth0em \kern-0.20em \rm Q}
        }}
\def\sCC{{\kern 0.27em \vrule height1.45ex width0.03em depth0em
          \kern-0.30em \rm C}}
\def\complex{{\mathchoice
    {\sCC}
    {\sCC}
  {\kern 0.225em \vrule height1.05ex width0.025em depth0em \kern-0.25em \rm C}
  {\kern 0.180em \vrule height0.78ex width0.020em depth0em \kern-0.20em \rm C}
        }}

\newcommand{\ba}{\begin{array}}
\newcommand{\ea}{\end{array}}
\newcommand{\be}{\begin{equation}}
\newcommand{\ee}{\end{equation}}
\newcommand{\bea}{\begin{eqnarray}}
\newcommand{\eea}{\end{eqnarray}}
\newcommand{\beaa}{\begin{eqnarray*}}
\newcommand{\eeaa}{\end{eqnarray*}}

\def\z{\zeta}

\def\qed{ \hfill \vrule width.25cm height.25cm depth0cm\smallskip}
\font\tenmath=msbm10 \font\sevenmath=msbm7 \font\fivemath=msbm5
\newfam\mathfam \textfont\mathfam=\tenmath
\scriptfont\mathfam=\sevenmath \scriptscriptfont\mathfam=\fivemath

\def \={{\buildrel {\rm (law)} \over =}}

\def\cB{\mathcal{B}}

\def\cD{\mathcal{D}}
\def\cE{\mathcal{E}}
\def\cF{\mathcal{F}}

\def\cH{\mathcal{H}}

\def\cP{\mathcal{P}}

\def\cS{\mathcal{S}}

\def\bR{\mathbb{R}}

\def\bE{\mathbb{E}}

\def\bC{\mathbb{C}}

\newcommand{\basa}{\begin{assumption}}
\newcommand{\easa}{\end{assumption}}

\newcommand{\bas}{\begin{assum}}
\newcommand{\eas}{\end{assum}}

\def\liminf{\mathop{\underline{\rm lim}}}

\def\bE{{\bf E}}

\newcommand{\ignore}[1]{}
\begin{document}

\renewcommand{\thefootnote}{\fnsymbol{footnote}}

\title[Stochastic Heat Equation with a Fractional-Colored Noise]
{The Stochastic Heat Equation with a Fractional-Colored Noise:
Existence of the Solution}

\author{Raluca M. Balan}

\address{Corresponding author. Department of Mathematics and Statistics,
University of Ottawa, 585 King Edward Avenue, Ottawa, Ontario K1N
6N5 Canada} \email{rbalan@uottawa.ca}

\author{Ciprian A. Tudor}

\address{SAMOS/MATISSE,
Centre d'Economie de La Sorbonne,\\ Universit\'e de
Panth\'eon-Sorbonne Paris 1, 90 rue de Tolbiac, 75634 Paris Cedex
13, France.} \email{tudor@univ-paris1.fr}



\thanks{{\em Acknowledgement of support:} The first author was supported by a grant from
the Natural Sciences and Engineering Research Council of Canada.}

\begin{abstract}
\noindent In this article we consider the stochastic heat equation
$u_{t}-\Delta u=\dot B$ in $(0,T) \times \bR^d$, with vanishing
initial conditions, driven by a Gaussian noise $\dot B$ which is
fractional in time, with Hurst index $H \in (1/2,1)$, and colored
in space, with spatial covariance given by a function $f$. Our
main result gives the necessary and sufficient condition on $H$
for the existence of the process solution. When $f$ is the Riesz
kernel of order $\alpha \in (0,d)$ this condition is
$H>(d-\alpha)/4$, which is a relaxation of the condition $H>d/4$
encountered when the noise $\dot B$ is white in space. When $f$ is
the Bessel kernel or the heat kernel, the condition remains
$H>d/4$.
\end{abstract}

\maketitle

\date{March 2, 2007}

\section{Introduction and Preliminaries}

Stochastic partial differential equations (s.p.d.e.'s) perturbed
by noise terms which bear a ``colored'' spatial covariance
structure (but remain white in time) have become increasingly
popular in the recent years, after the fundamental work of
\cite{Da1}. Such an equation can be viewed as a more flexible
alternative to a classical s.p.d.e. driven by a space-time white
noise, and therefore it can be used to model a more complex
physical phenomenon which is subject to random perturbations. The
major drawback of this theory is that it is mathematically more
challenging than the classical theory, usually relying on
techniques from potential analysis. One advantage is that an
s.p.d.e. perturbed by a colored noise possesses a process solution
(under relatively mild conditions on the covariance structure), in
contrast with its white-noise driven counterpart, for which the
solution is well understood only in the sense of distributions.
Another advantage is the fact that such an equation can lead to a
better understanding of a complex physical situation.

This article continues the line of research initiated by
\cite{Da1}, the focus being on a relatively simple s.p.d.e., the
stochastic heat equation. The novelty comes from the fact that
random noise perturbing the equation possesses a colored temporal
structure given by the covariance of a fractional Brownian motion
(fBm), along with the colored spatial structure of \cite{Da1}.

We recall that a fBm on the real line is a zero-mean Gaussian
process with covariance function
$R_{H}(t,s)=(t^{2H}+s^{2H}-|t-s|^{2H})/2$, where $H \in (0,1)$.
There is a huge amount of literature dedicated to the fBm, due to
its mathematical tractability, and its many applications. We refer
the reader to \cite{N2} for a comprehensive review on this
subject.

Recently, the fBM made its entrance in the area of s.p.d.e.'s. We
refer, among others, to \cite{MaNu}, \cite{TTV}, \cite{DMD},
\cite{NuOuk}, \cite{gubinelli-lejay-tindel06} or
\cite{quersardanyons-tindel06}. For example, in the case of the
stochastic heat equation driven by a Gaussian noise which is
fractional  in time and with a rather general covariance in space,
it has been proved in \cite{TTV} that if the time variable belongs
to $[0,T]$ and the space variable belongs to $S^{1}$ (the unit
circle) then the solution exists if and only if $H>1/4$. The case
of the same equation driven by a fractional-white noise with space
variable in $\mathbb{R}^{d}$ has been treated in \cite{MaNu}, and
it follows that a process solution exists if $H>d/4$.


We note in passing that very few results are available in the
literature, in the case of non-linear equations driven by a
fractional Gaussian noise, due to difficulties encountered in the
stochastic calculus associated to this noise. To circumvent these
difficulties, a new pathwise method has been developed recently in
\cite{gubinelli-lejay-tindel06} and
\cite{quersardanyons-tindel06}, treating the nonlinear stochastic
heat equation, respectively the nonlinear stochastic wave
equation. In both these articles, the noise term carries a colored
spatial covariance structure as well, which is the situation that
we investigate in the present paper too.

In the present article we consider the (linear) stochastic heat
equation in the domain $[0,T] \times \bR^d$ driven by a Gaussian
noise $B$ which has a fractional time component, of Hurst index $H
\in (1/2,1)$, and a colored spatial component. Therefore, our work
lies at the intersection of the two different lines of research
mentioned above, namely those developed in \cite{Da1},
respectively \cite{NuOuk}. The solution of the equation will be
given in the mild formulation, but can also be viewed as a
distribution solution. Therefore, the first step we need to take
is to develop a stochastic calculus with respect to the noise $B$.
Since our equation is linear, only spaces of deterministic
integrands are considered in the present paper. Our main result
identifies the necessary and sufficient condition on the Hurst
index $H$ for the existence of the process solution. We should
mention that in this case, the solution is a Gaussian process, and
hence Reproducing Kernel Hilbert Space techniques can be used to
investigate its properties. This will be the subject of future
work.

In preparation for treating the existence problem in its full
generality, we studied first the case of the heat equation
perturbed by a Gaussian noise which is fractional in time, but
white in space (Section \ref{sect-fract-white}). In this case, it
is known from \cite{MaNu} that the process solution exists if
$H>d/4$, which forces a spatial dimension $d\in \{1,2,3\}$; in the
present article, we strengthen this result by proving that the
condition $H>d/4$ is in fact necessary for the existence of the
solution. In contrast, when a spatial covariance structure is
embedded in the noise (Section \ref{sect-fract-color}) the
condition for the existence of the solution can be relaxed so that
it imposes no restrictions on the spatial dimension $d$. When the
color in space is given by the Riesz kernel of order $\alpha$, we
prove that the necessary and sufficient condition for the
existence of the solution is $H>(d-\alpha)/4$. This demonstrates
that a suitable choice of the spatial covariance structure can
compensate for the drawbacks of the fractional component. However,
it turns out that if the spatial covariance is given by the Bessel
or the heat kernel, the condition remains $H>d/4$, whereas for the
Poisson kernel the condition becomes $H>(d+1)/4$.

This article contains 3 appendices. Appendix A contains a lemma
which is heavily used in the present paper. This lemma is the tool
which allows us to import the Fourier transform techniques from
$\bR$ to the bounded domain $[0,T]$. (So far, this type of
techniques have been exploited only on $\bR$; see e.g.
\cite{PT1}.) Appendix B contains the proof of a technical
statement. Appendix C contains a result which is essentially due
to \cite{Da1}; we include it since we could not find a direct
reference.

\vspace{3mm}

We begin now to introduce the notation that will be used
throughout this paper.

If $U \subset \bR^n$ is an open set, we denote by $\cD(U)$ the
space of all infinitely differentiable functions whose support is
compact and contained in $U$. By $\cD'(U)$ we denote the set of
continuous linear functionals on $\cD(U)$ which is known as the
space of {\it distributions}. We let $\cS(\bR^n)$ be the Schwarz
space of all decreasing functions on $\bR^n$ and $\cS'(\bR^n)$ be
the space of {\it tempered distributions}, i.e. continuous linear
functionals on $\cS(\bR^n)$. For an arbitrary function $g$ on
$\bR^d$ the translation by $x$ is denoted by $g_{x}$, i.e.
$g_{x}(y)=g(x+y)$. The reflection by zero is denoted by
$\tilde{g}$, i.e. $\tilde{g}(x)= g(-x)$.

 For any function $\phi \in \cS(\bR^n)$
we define its Fourier transform by
\begin{equation*}
\cF\phi (\xi) = \int_{\bR^d} \exp ( -i \xi \cdot x ) \phi (x) dx.
\end{equation*}
The map $\cF:\cS(\bR^n) \to \cS(\bR^n) $ is an isomorphism which
extends uniquely to a unitary isomorphism of $L_{2}(\bR^n)$; this
map can also be extended to $\cS'(\bR^n)$. We define the
convolution $(f*g)(x)= \linebreak \int_{R^n}f(x-y)g(y)dy$.

For any interval $(a,b)\subset \bR$ and $\varphi \in L_2(a,b)$, we
define the {\em restricted Fourier transform} of $\varphi$ with
respect to $(a,b)$ by:
$${\cF}_{a,b} \varphi (\tau)=\int_{a}^{b}e^{-i\tau t}\varphi(t)dt.$$

As in \cite{PT}, if $f \in L_{1}(0,T)$ and $\alpha>0$, we define
the {\em fractional integral} of $f$ of order $\alpha$ by:
$$(I_{T-}^{\alpha}f)(t)=\frac{1}{\Gamma(\alpha)}
\int_{t}^{T}(u-t)^{\alpha-1} f(u)du.$$

Finally, we denote by $\cB_{b}(\bR^d)$ the class of all bounded
Borel sets in $\bR^d$.

\section{The Fractional-White Noise}
\label{sect-fract-white}

The purpose of this section is to identify the necessary and
sufficient condition for the existence of a distribution-solution
of the stochastic heat equation, driven by a Gaussian noise which
is {\em fractional in time} and {\em white in space}.

 Our main
result is similar to Theorem 11, \cite{Da1}, which identifies the
necessary and sufficient condition for the existence of a
distribution-solution of an arbitrary s.p.d.e.'s driven by a
Gaussian noise which is {\em white in time} and {\em colored in
space}.

 In the first subsection we examine some spaces of
deterministic integrands, which are relevant for the stochastic
calculus with respect to fractional processes, and we explore the
connection with Dalang's theory of stochastic integration. In the
second subsection, we describe the Gaussian noise and its
stochastic integral. In the third subsection, we introduce the
process-solution and the distribution-solution of the stochastic
heat equation driven by this noise, and we identify the necessary
and sufficient condition for the existence of these processes.

\subsection{Spaces of deterministic integrands}

We begin by introducing the usual spaces associated with the
fractional temporal noise. Throughout this article we suppose that
$H \in (1/2,1)$ and we let $\alpha_{H}=H(2H-1)$.

Let $\cH(0,T)$ be the completion of $\cD(0,T)$ with respect to the
inner product
\begin{eqnarray*}
\langle \varphi,\psi \rangle_{\cH(0,T)} & = & \alpha_{H}
\int_{0}^{T} \int_{0}^{T}\varphi(u) |u-v|^{2H-2} \psi(v)dvdu \\
& = & \alpha_H c_{H} \int_{\bR} \cF_{0,T}\varphi(\tau)
\overline{\cF_{0,T} \psi(\tau)} |\tau|^{-(2H-1)}d\tau,
\end{eqnarray*}
where the second equality follows by Lemma \ref{LemmaA1}.(b) with
$c_{H}=[2^{2(1-H)}\pi^{1/2}]^{-1} \Gamma(H-1/2)/\Gamma(1-H)$. Note
that $c_{H}=q_{2H-1}$, where the constant $q_{\alpha}$ is defined
in Lemma \ref{LemmaA1} (Appendix A).

Let $\cE(0,T)$ be the class of all linear combinations of
indicator functions $1_{[0,t]}, t \in [0,T]$. One can see that
$\cH(0,T)$ is also the completion of $\cE(0,T)$ with respect to
the inner product $$\langle 1_{[0,t]},
1_{[0,s]}\rangle_{\cH(0,T)}= R_H(t,s).$$ (To see this, note that
every $1_{[0,a]} \in \cE(0,T)$ there exists a sequence $(\varphi
_{n})_{n} \subset \cD(0,T)$ such that $\varphi_n(t)\to
1_{[0,a]}(t), \forall t$ and ${\rm supp} \ \varphi_n \subset K$
for all $n$, where $K \subset (0,T)$ is a compact. By the
dominated convergence theorem, it follows that $\|
\varphi_{n}-1_{[0,t]} \|_{\cH(0,T)}\to 0$.)

\vskip0.5cm

In is important to emphasize that the space $\cH(0,T)$ may contain
{\em distributions}. We will justify this statement, using the
argument of \cite{Da1}. First, note that $\cH(0,T) \subset
\cH(\bR)$, where $\cH(\bR) $ is the completion of $\cD(\bR) $ with
respect to the inner product
\begin{equation*}
\langle \varphi, \psi \rangle_{\cH(\bR)}= \alpha_{H} \int _{\bR}
\int_{\bR} \varphi (u) \psi (v) |u-v|^{2H-2} dvdu.
\end{equation*}
The space $\cH(\bR)$ appears in several papers treating colored
noises. (In fact $\cH(\bR)$ is a particular instance of the space
$\cP^{(d)}_{0,x}$ of \cite{Balan}, in the case $\mu (d\xi
)=|\xi|^{-(2H-1)}d\xi$ and $d=1$.) From p. 9 of \cite{Da1}, we
know that
\begin{equation}
\label{cH-included-barcH} \cH(\bR) \subset
\overline{\cH}(\bR):=\{S \in {\cS}'(\bR); {\cF} S \ \mbox{is a
function} , \ \int_{\bR} |\cF S(\tau)|^{2}|\tau|^{-(2H-1)}d \tau<
\infty\}.
\end{equation}
 Since $|\tau|^{2}<1+|\tau|^{2}$, one can
easily see that $\overline{\cH}(\bR) \subset
{\cH}^{-(H-1/2)}(\bR)$, where
$${\cH}^{-(H-1/2)}(\bR) :=  \{S \in \cS'(\bR); {\cF}S
\ \mbox{is a function}, \ \int_{\bR}|\cF
S(\tau)|^{2}(1+|\tau|^{2})^{-(H-1/2)}d\tau< \infty\}$$ is the
fractional Sobolev space of index $-(H-1/2)$ (see p. 191,
\cite{folland95}). Therefore, the elements of $\cH(\bR)$ are
tempered distributions on $\bR$ of negative order $-(H-1/2)$.
(This was also noticed by several authors; see e.g. p. 9,
\cite{N2}, or \cite{PT1}, \cite{PT})

On the other hand, similarly to (\ref{cH-included-barcH}), one can
show that
\begin{equation}
\label{H(0,T)-included-barH(0,T)} \cH(0,T) \subset
\overline{\cH}(0,T):= \{S \in {\cS}'(\bR); {\cF}_{0,T} S \
\mbox{is a function} , \ \int_{\bR} |\cF_{0,T}
S(\tau)|^{2}|\tau|^{-(2H-1)}d \tau< \infty\},
\end{equation}
 where {\em the restricted Fourier transform}
${\cF}_{0,T}S$ of a distribution $S \in \cS'(\bR)$ is defined by:
$\langle {\cF}_{0,T}S, \varphi\rangle =\langle
S,{\cF}_{0,T}\varphi \rangle, \ \forall \varphi \in \cS(\bR)$.

\begin{remark}
{\rm  One can show that (see p.10, \cite{N2})
\begin{equation}
\label{L2(0,T)-subset-H(0,T)} L_{2}(0,T) \subset L_{1/H}(0,T)
\subset |\cH(0,T)| \subset \cH(0,T),
\end{equation}
where $|\cH(0,T)|=\{f:[0,T]\times \mathbb{R}^{d} \to \bR \
\mbox{measurable}; \ \int_{0}^{T} \int_{0}^{T} |f(u)| |f(v)| \vert
\vert u-v\vert ^{2H-2} dudv<\infty\}$. }
\end{remark}

A different approach of characterizing the space $\cH(0,T)$ is
based on the {\em transfer operator}. We recall that the kernel
$K_{H}(t,s),t>s$ of the fractional covariance function $R_{H}$ is
defined by:
$$K_{H}(t,s)=c_{H}^*s^{1/2-H}
\int_{s}^{t}(u-s)^{H-3/2}u^{H-1/2}du,$$ where $c_{H}^*=\{\alpha_H
\Gamma(3/2-H)/[\Gamma(2-2H) \Gamma(H-1/2)]\}^{1/2}$.

 Note that
$R_{H}(t,s)=\int_{0}^{t \wedge s}K_{H}(t,u)K(s,u)du$ (see p.7-8,
\cite{N2}) and hence
$$\langle K_{H}^* 1_{[0,t]}, K_{H}^*1_{[0,s]}
\rangle_{L_{2}(0,T)}= \int_{0}^{t \wedge s} K_{H}(t,u)K_{H}(s,u)du
=R_H(t,s) =\langle 1_{[0,t]},1_{[0,s]}\rangle_{\cH(0,T)},$$
 i.e. $K_{H}^*$ is an
isometry between $(\cE(0,T), \langle \cdot,\cdot
\rangle_{\cH(0,T)})$ and $L_{2}(0,T)$. Since $\cH(0,T)$ is the
completion of $\cE(0,T)$ with respect to $\langle \cdot, \cdot
\rangle_{\cH(0,T)}$, this isometry can be extended to $\cH(0,T)$.
We denote this extension by $K_{H}^*$. In fact, one can prove that
the map $K_{H}^* : \cH(0,T) \rightarrow L_{2}(0,T)$ is also
surjective (see the proof of Lemma \ref{isometry-cH-L2}).

\vspace{3mm}

We are now introducing the space of deterministic integrands
associated with the a noise which is fractional in time and white
in space. This space was also considered in \cite{NuOuk}.

More precisely, let $\cH$ be the completion of $\cD((0,T) \times
\bR^d)$ with respect to the inner product
\begin{eqnarray*}
\langle \varphi,\psi \rangle_{\cH} & = & \alpha_{H} \int_{0}^{T}
\int_{0}^{T} \int_{\bR^d}\varphi(u,x) |u-v|^{2H-2} \psi(v,x)dx
dv du\\
& = & \alpha_H c_{H} \int_{\bR}\int_{\bR^d}
\cF_{0,T}\varphi(\tau,x) \overline{\cF_{0,T} \psi(\tau,x)}
|\tau|^{-(2H-1)}dx d\tau  \\
& = & \int_{\bR^d} \langle \varphi(\cdot,x), \psi(\cdot,x)
\rangle_{\cH(0,T)}dx
\end{eqnarray*}
where the second equality above follows by Lemma
\ref{LemmaA1}.(b), and the third is due to Fubini's theorem.

If we let $\cE$ be the space of all linear combinations of
indicator functions $1_{[0,t] \times A}, t \in [0,T], A \in
\cB_{b}(\bR^d)$, then one can prove that $\cH$ is also the
completion of $\cE$ with respect to the inner product
$$\langle 1_{[0,t] \times A}, 1_{[0,s] \times B}\rangle_{\cH}=R_H(t,s)
 \lambda(A \cap B)$$ where $\lambda$ is the Lebesgue measure
on $\bR^d$. (The argument is similar to the one used in the
temporal case.) Similarly to (\ref{H(0,T)-included-barH(0,T)}),
one can show that
$$\cH \subset \overline{\cH}:=\{S:\bR^d \to {\cS}'(\bR);
 {\cF}_{0,T} S(\cdot,x) \ \mbox{is a
function} \ \forall x \in \bR^d, \ (\tau,x) \mapsto
\cF_{0,T}S(\tau,x) \ \mbox{is measurable},$$ $$ \ \mbox{and} \
\int_{\bR} \int_{\bR^d} |\cF_{0,T}
S(\tau,x)|^{2}|\tau|^{-(2H-1)}dx d \tau < \infty\}.$$

\noindent Using (\ref{L2(0,T)-subset-H(0,T)}) and the fact that
$$\|S \|_{\cH}^{2}=\int_{\bR^d} \|S(\cdot,x) \|_{\cH(0,T)}^{2}dx,
\quad \forall S \in \cH,$$ one can show that
$$L_{2}((0,T) \times \bR^d) \subset |\cH| \subset \cH \subset L_{2}(\bR^d;
\cH(0,T)),$$ where $|\cH|=\{\varphi:[0,T]\times \bR^d \
\mbox{measurable}; \  \int_{0}^{T} \int_{0}^{T} \int_{\bR^d}
|\varphi(u,x)| |\varphi(v,x)| |u-v|^{2H-2} dx dv du<\infty\}$.

\vskip0.5cm

The next result gives an alternative criterion for verifying that
a function $\varphi$ lies in $\cH$. (It can be compared to
Theorems 2 and 3 of \cite{Da1}.)

\begin{theorem}
\label{criterion-cH} Let $\varphi :[0,T] \times \bR^d \rightarrow
\bR$ be a function which satisfies the following conditions:

(i) $\varphi(\cdot,x) \in L_{2}(0,T)$ for every $x \in \bR^d$;

(ii) $(\tau,x) \mapsto \cF_{0,T} \varphi (\tau,x)$ is measurable;

(iii) $\int_{\bR} \int_{\bR^d}  |\cF_{0,T} \varphi(\tau,x)|^{2}
|\tau|^{-(2H-1)} dx d\tau<\infty$.

\noindent Then $\varphi \in \cH$.
\end{theorem}

\begin{proof}
Similarly to Proposition 3.3, \cite{PT1}, we let
\begin{eqnarray*}
\tilde \Lambda& =& \{\varphi:[0,T] \times \bR^d \to \bR;
\varphi(\cdot,x) \in L_{2}(0,T) \ \forall x, (\tau,x) \mapsto
\cF_{0,T} \varphi (\tau,x) \ \mbox{is measurable, and} \\
& & \|\varphi \|_{\tilde \Lambda}^{2}:=c_{1} \int_{\bR}
\int_{\bR^d}|\cF_{0,T} \varphi(\tau,x)|^{2} |\tau|^{-(2H-1)} dx
d\tau<\infty\}
\\
\Lambda&=&\{\varphi:[0,T] \times \bR^d \to \bR; \ \|\varphi
\|_{\Lambda}^{2}:=c_{2}\int_{0}^{T} \int_{\bR^d}
[I_{T-}^{H-1/2}(u^{H-1/2}\varphi(u,x))(s)]^{2} s^{-(2H-1)} dx ds
<\infty\}
\end{eqnarray*}
where $c_{1}=\alpha_Hc_{H}$ and
$c_2=\{c_{H}^*\Gamma(H-1/2)\}^{2}$. We now prove that
\begin{equation}
\label{tildeLambda-in-Lambda} \tilde \Lambda \subset \Lambda \quad
\mbox{and} \quad \|\varphi \|_{\tilde \Lambda}=\|\varphi
\|_{\Lambda}, \quad \forall \varphi \in \tilde \Lambda.
\end{equation}

Let $\varphi \in \tilde \Lambda$ be arbitrary. Since
$\varphi(\cdot,x) \in L_{2}(0,T) \subset \cH(0,T)$  and
$K_{H}^{*}$ is an isometry from $\cH(0,T)$ to $L_{2}(0,T)$, we
have $\|\varphi(\cdot,x)\|_{\cH(0,T)}^{2}=\|K_{H}^{*}\varphi
(\cdot,x)\|_{L_{2}(0,T)}^{2}$ for all $x \in \bR^d$, that is
$$c_{1}\int_{\bR}
|\cF_{0,T} \varphi(\tau,x)|^{2} |\tau|^{-(2H-1)} d\tau
=c_{2}\int_{0}^{T} [I_{T-}^{H-1/2}(u^{H-1/2}\varphi(u,x))(s)]^{2}
s^{-(2H-1)} ds, \quad \forall x \in \bR^d.$$ Integrating with
respect to $dx$ and using Fubini's theorem, we get $\|\varphi
\|_{\tilde \Lambda}=\|\varphi \|_{\Lambda}$. The fact that
$\|\varphi \|_{\tilde \Lambda}<\infty$ forces $\|\varphi
\|_{\Lambda}<\infty$, i.e. $\varphi \in \Lambda$. This concludes
the proof of (\ref{tildeLambda-in-Lambda}).

Next we prove that
\begin{equation}
\label{E-dense-in-Lambda} \cE \ \mbox{is dense in} \ \Lambda \
\mbox{with respect to} \ \|\cdot \|_{\Lambda}.
\end{equation}

Let $\varphi \in \Lambda$ and $\varepsilon>0$ be arbitrary. Since
the map $(s,x) \mapsto I_{T-}^{H-1/2}(u^{H-1/2}\varphi(u,x))(s)$
belongs to $L_{2}((0,T) \times \bR^d,d\lambda_H \times dx)$ where
$\lambda_{H}(s)=s^{-(2H-1)}ds$, there exists a simple function
$g(s,x)=\sum_{k=1}^{n}b_{k}1_{[c_k,d_k)}(s)1_{A_k}(x)$ on $(0,T)
\times \bR^d$, with $b_{k} \in \bR$, $0<c_k<d_k<T$ and $A_{k}
\subset \bR^d$ Borel set, such that
\begin{equation}
\label{first-approx-IT-varphi}
\int_{0}^{T}\int_{\bR^d}[I_{T-}^{H-1/2}(u^{H-1/2}\varphi(u,x))(s)-g(s,x)]^{2}s^{-(2H-1)}dxds
<\varepsilon.
\end{equation}

\noindent By relation (8.1) of \cite{PT}, we know that there
exists an elementary function $l_{k} \in \cE(0,T)$ such that
$$\int_{0}^{T}
[1_{[c_k,d_k)}(s)-I_{T-}^{H-1/2}(u^{H-1/2}l_k(u))(s)]^{2}s^{-(2H-1)}ds<\varepsilon/C_g,$$
where we chose $C_g:=n\sum_{k=1}^{n}b_{k}^{2}\lambda(A_k)$. Define
$l(s,x)=\sum_{k=1}^{n}b_k l_{k}(t)1_{A_k}(x)$ and note that $l \in
\cE$. Then
\begin{equation}
\label{second-approx-IT-varphi} \int_{0}^{T} \int_{\bR^d}
[g(s,x)-I_{T-}^{H-1/2}(u^{H-1/2}l(u,x))(s)]^{2}s^{-(2H-1)}ds<\varepsilon.
\end{equation}

\noindent From (\ref{first-approx-IT-varphi}) and
(\ref{second-approx-IT-varphi}), we get
$$\int_{0}^{T}\int_{\bR^d}[I_{T-}^{H-1/2}(u^{H-1/2}\varphi(u,x))(s)-
I_{T-}^{H-1/2}(u^{H-1/2}l(u,x))(s)]^{2}s^{-(2H-1)}dxds
<4\varepsilon,$$ i.e. $\|\varphi-l \|_{\Lambda}^{2}<4 \varepsilon
c_2$. This concludes the proof of (\ref{E-dense-in-Lambda}).

From (\ref{tildeLambda-in-Lambda}) and (\ref{E-dense-in-Lambda}),
we infer immediately that $\cE$ is dense in $\tilde \Lambda$ with
respect to $\|\cdot \|_{\tilde \Lambda}$. Since $\|\cdot
\|_{\tilde \Lambda}=\|\cdot \|_{\cH}$ and $\cH$ is the completion
of $\cE$ with respect to $\|\cdot \|_{\cH}$, it follows that
$\tilde \Lambda \subset \cH$. This concludes the proof of the
theorem. \qed
\end{proof}

\vspace{3mm}

 As in \cite{NuOuk}, we define the transfer operator by:
\begin{equation}
\label{def-K-H*} (K_{H}^*1_{[0,t] \times
A})(s,x):=K_{H}(t,s)1_{[0,t] \times A}(s,x)
.
\end{equation}

\noindent Note that
\begin{eqnarray*}
\langle K_{H}^* 1_{[0,t] \times A}, K_{H}^*1_{[0,s] \times B}
\rangle_{L_{2}((0,T) \times \bR^d)}& =& \left( \int_{0}^{t \wedge
s} K_{H}(t,u)K_{H}(s,u)du \right) \langle 1_A, 1_B
\rangle_{L_{2}(\bR^d)} \\
& = & R_H(t,s)\lambda(A \cap B) =\langle 1_{[0,t] \times
A},1_{[0,s] \times B}\rangle_{\cH},
\end{eqnarray*}
i.e. $K_{H}^*$ is an isometry between $(\cE, \langle \cdot,\cdot
\rangle_{\cH})$ and $L_{2}((0,T) \times \bR^d)$. Since $\cH$ is
the completion of $\cE$ with respect to $\langle \cdot, \cdot
\rangle_{\cH}$, this isometry can be extended to $\cH$. We denote
this extension by $K_{\cH}^*$.

\begin{lemma}
\label{isometry-cH-L2} $K_{\cH}^* : \cH \rightarrow L_{2}((0,T)
\times \bR^d)$ is surjective.
\end{lemma}

\begin{proof}
Note that $1_{[0,t] \times A} \in K_{\cH}^{*}(\cH)$ for all $t \in
[0,T], A \in \cB_{b}(\bR^d)$; hence $\varphi \in K_{\cH}^{*}(\cH)$
for every $\varphi \in \cE$. Let $f \in L_{2}((0,T) \times \bR^d)$
be arbitrary. Since $\cE$ is dense in $L_{2}((0,T) \times \bR^d)$,
there exists a sequence $(f_{n})_{n} \subset \cE$ such that
$\|f_n-f \|_{L_{2}((0,T) \times \bR^d)} \to 0$. Since $f_n \in
K_{\cH}^{*}(\cH)$, there exists $\varphi_n \in \cH$ such that
$f_{n}=K_{\cH}^* \varphi_n$. The sequence $(\varphi_n)_n$ is
Cauchy in $\cH$:
$$\|\varphi_n-\varphi_m \|_{\cH}=\|K_{\cH}^* \varphi_n-K_{\cH}^*
\varphi_m \|_{L_{2}((0,T)\times \bR^d)}=\|f_n-f_m \|_{L_{2}((0,T)
\times \bR^d)} \rightarrow 0$$ as $m,n \to\infty$. Since $\cH$ is
complete, there exists $\varphi \in \cH$ such that $\|\varphi_n
-\varphi\|_{\cH} \to 0$. Hence
$\|f_{n}-K_{\cH}^*\varphi\|_{L_{2}((0,T) \times
\bR^d)}=\|K_{\cH}^*\varphi_n- K_{\cH}^*\varphi\|_{L_{2}((0,T)
\times \bR^d)} \to 0$. But $\|f_n-f \|_{L_{2}((0,T) \times \bR^d)}
\to 0$. We conclude that $K_{\cH}^*\varphi=f$. \qed
\end{proof}

\vspace{3mm}

\begin{remark}
{\rm Note that for every $\varphi \in \cE$
\begin{eqnarray}
(K_{H}^* \varphi)(s,x)&=&\int_s^T\varphi(r,x)\frac{\partial
K_H}{\partial r} (r,s)dr = c_{H}^*\int_s^T\varphi(r,x)
\left(\frac{r}{s}
\right)^{H-1/2}(r-s)^{H-3/2}dr \nonumber \\
\label{K-H*-as-fract-integr}
&=&c_{H}^*\Gamma\left(H-\frac{1}{2}\right) s^{-(H-1/2)}
I_{T-}^{H-1/2}(u^{H-1/2}\varphi(u,x))(s).
\end{eqnarray}
Using Lemma \ref{isometry-cH-L2}, we can {\em formally} say that
$$\cH=\{\varphi  \ \mbox{such that} \ (s,x) \mapsto
s^{-(H-1/2)} I_{T-}^{H-1/2}(u^{H-1/2}\varphi(u,x))(s)  \ \mbox{is
in} \ L_{2}((0,T) \times \bR^d) \}.$$ }
\end{remark}

\subsection{The Noise and the Stochastic Integral}
\label{fractionar-white-noise}

In this paragraph we describe the Gaussian noise which is randomly
perturbing the heat equation. This noise is assumed to be
fractional in time and white in space and was also considered by
other authors (see \cite{NuOuk}).

Let $F=\{F(\varphi); \varphi \in {\cD}((0,T) \times \bR^{d})\}$ be
a zero-mean Gaussian process with covariance
\begin{equation}
\label{noise1}E(F(\varphi)F(\psi))= \langle \varphi, \psi
\rangle_{\cH}.
\end{equation}
Let $H^F$ be the Gaussian space of $F$, i.e. the closed linear
span of $\{F(\varphi);\varphi \in {\cD}((0,T) \times \bR^{d})\}$
in $L^{2}(\Omega)$.

For every indicator function $1_{[0,t] \times A} \in \cE$, there
exists a sequence $(\varphi_n)_{n} \subset \cD((0,T) \times
\bR^d)$ such that $\varphi_n \to 1_{[0,t] \times A}$ and ${\rm
supp} \ \varphi_n \subset K, \ \forall n$, where $K \subset (0,T)
\times \bR^d$ is a compact. Hence $\| \varphi_n-1_{[0,t] \times A}
\|_{\cH} \to 0$ and $E(F(\varphi_{m})-F(\varphi_{n}))^2= \|
\varphi_{m}-\varphi_{n} \|_{\cH} \to 0$ as $m,n \rightarrow
\infty$, i.e. the sequence $\{F(\varphi_{n})\}_{n}$ is Cauchy in
$L_{2}(\Omega)$. A standard argument shows that its limit does not
depend on $\{\varphi_{n}\}_{n}$. We set $F_{t}(A)=F(1_{[0,t]
\times A})=_{L_{2}(\Omega)}\lim_{n}F(\varphi_{n}) \in H^{F}$.
We extend $F$ by linearity to $\cE$. A limiting argument and
relation (\ref{noise1}) shows that
$$E(F(\varphi)F(\psi))=\langle \varphi,\psi \rangle_{\cH}, \quad
\forall \varphi,\psi \in {\cE},$$ i.e. $\varphi \mapsto
F(\varphi)$ is an isometry between $(\cE, \langle \cdot, \cdot
\rangle_{\cH})$ and $H^F$. Since $\cH$ is the completion of $\cE$
with respect to $\langle \cdot, \cdot \rangle_{\cH}$, this
isometry can be extended to $\cH$, giving us the stochastic
integral with respect to $F$. We will use the notation
$$F(\varphi)=\int_{0}^{T} \int_{\bR^d}\varphi(t,x)F(dt,dx).$$

\begin{remark}
{\rm One can use the transfer operator $K_{\cH}^*$ to explore the
relationship between $F(\varphi)$ and Walsh's stochastic integral
(introduced in \cite{walsh86}). More precisely, using Lemma
\ref{isometry-cH-L2}, we define
\begin{equation}
\label{def-of-W} W(\phi):=F((K_{\cH}^*)^{-1}(\phi)),  \ \ \ \phi
\in L_{2}((0,T) \times \bR^d).
\end{equation} Note that
$$\bE(W(\phi)W(\eta))=\langle (K_{\cH}^*)^{-1}(\phi),
(K_{\cH}^*)^{-1}(\eta)\rangle_{\cH}=\langle \phi,\eta
\rangle_{L_{2}((0,T) \times \bR^d)},$$ i.e. $W=\{W(\phi); \phi \in
L_2((0,T) \times \bR^{d})\}$ is a space-time white noise. Using
the stochastic integral notation, we write $W(\phi)=\int_0^T
\int_{\bR^d} \phi(t,x)W(dt,dx)$ for all $\phi \in L_{2}((0,T)
\times \bR^d)$. (Note that $W(\phi)$ is Walsh's stochastic
integral with respect to the noise $W$.)

Let $H^W$ be the Gaussian space of $W$, i.e. the closed linear
span of $\{W(\phi); \phi \in L_2((0,T) \times \bR^{d})\}$ in
$L_{2}(\Omega)$. By (\ref{def-of-W}), we can see that
$H^{W}=H^{F}$. The following diagram summarizes these facts:

\begin{picture}(200,120)(10,10)
\put(20,90){\vector(0,-1){70}} \put(15,95){$\cH$}
\put(110,90){\vector(-1,-1){70}} \put(105,95){$L_{2}((0,T) \times
\bR^d)$} \put(30,100){\vector(1,0){70}} \put(10,10){$H^{F}=H^{W}$}
\put(50,105){$K_{\cH}^{*}$} \put(10,50){$F$} \put(80,50){$W$}
\put(140,70){$F(\varphi)=W(K_{\cH}^*\varphi), \ \ \ \forall
\varphi \in \cH, \ {\rm i.e.}$}
 \put(140,40){$\int_{0}^{T}\int_{\bR^d}\varphi(t,x)F(dt,dx)=\int_{0}^{T}
\int_{\bR^d} (K_{\cH}^*\varphi)
 (t,x)W(dt,dx),\ \ \ \forall \varphi \in \cH.$}
\end{picture}

\noindent In particular, $F(t,A)=\int_{0}^{t}
\int_{A}K_{H}(t,s)W(ds,dy)$. This relationship will not be used in
the present paper.

}
\end{remark}

\subsection{The Solution of the Stochastic Heat Equation}
\label{fractionar-white-sol}

We consider the stochastic heat equation driven by the noise $F$,
written {\em formally} as:

\begin{equation}
\label{heat-eq-H}v_{t}-\Delta v =  \dot F , \ \ \ \mbox{in} \
(0,T) \times \mathbb{R}^{d}, \ \ \  v(0,\cdot)= 0,
\end{equation}
where $\Delta v$ denotes the Laplacian of $v$, and $v_{t}$ is the
partial derivative with respect to $t$.

Let $G$ be the fundamental solution of the classical heat
equation, i.e.
\begin{equation}
\label{fund-sol-heat} G(t,x)=\left\{
\begin{array}{ll} (4 \pi t)^{-d/2} \exp\left( -\frac{|x|^{2}}{4t}\right) & \mbox{if $t>0, x \in R^{d}$} \\
0 & \mbox{if $t \leq 0, x \in R^{d}$}
\end{array} \right.
\end{equation}
 Let
$g_{tx}(s,y):=(G_{tx})^{\verb2~2}(s,y)=G(t-s, x-y)$. The following
result shows that the kernel $G$ has the desired regularity, which
is needed to apply Theorem \ref{criterion-cH}.

\begin{lemma}
\label{eta*G-in-L2} If $\varphi=\eta * \tilde G$, where $\eta \in
\cD(0,T) \times \bR^d)$, then $$\varphi(\cdot,x) \in L_{2}(0,T)
\quad \forall x \in \bR^d.$$
\end{lemma}

\begin{proof}
Without loss of generality, we suppose that
$\eta(t,x)=\phi(t)\psi(x)$, where $\phi \in \cD(0,T)$ and $\psi
\in \cD(\bR^d)$. Using Minkowski's inequality for integrals (see
p. 271, \cite{stein70}), we have
\begin{eqnarray*}
\left(\int_{0}^{T}|\varphi(t,x)|^{2} dt \right)^{1/2}&=&
\left(\int_{0}^{T} \left| \int_{0}^{T}\int_{\bR^d}
\phi(s)\psi(y)G(s-t,y-x)1_{\{s>t\}}dyds \right|^{2} dt
\right)^{1/2} \\
& \leq &\int_{0}^{T}\int_{\bR^d}
|\phi(s)\psi(y)|\left(\int_{0}^{s} |G(s-t,y-x)|^{2}dt
\right)^{1/2}dy ds.
 \end{eqnarray*}

\noindent Using the change of variables $s-t=t'$ and $1/t'=u$, we
get
$$\int_{0}^{s}
|G(s-t,y-x)|^{2}dt=
\frac{1}{(4 \pi)^d}\int_{1/s}^{\infty}u^{d-2} e^{-|y-x|^2 u/2} du
\leq \frac{C}{|y-x|^{2(d-1)}},$$ and
$$\left(\int_{0}^{T}|\varphi(t,x)|^{2} dt \right)^{1/2} \leq
C \int_{0}^{T} |\phi(s)|
\int_{\bR^d}\frac{|\psi(y)|}{|y-x|^{d-1}}dyds \leq C (|\psi|*
R_1)(x)<\infty,$$ where $R_{1}(x)=\gamma_{1,d}|x|^{-(d-1)}$ is the
Riesz kernel of order $1$ in $\bR^d$ and the convolution $|\psi|
* R_{1}$ is well-defined by Theorem V.1, p. 119, \cite{stein70}.
\qed
\end{proof}

\vspace{3mm}

The next result gives the necessary and sufficient condition for
the existence of the process solution.

\begin{theorem}
\label{H-psi*Gtilde} If
\begin{equation}
\label{cond-H} H > \frac{d}{4} \ ,
\end{equation}
then: (a) $g_{tx} \in |\cH|$ for every $(t,x) \in [0,T] \times
\mathbb{R}^{d}$; (b) $\eta * \tilde G \in \cH$ for every $\eta \in
\cD((0,T) \times \bR^d)$. Moreover,
\begin{equation}
\label{nec-suf-cond-H} \|g_{tx}\|_{\cH} <\infty \ \forall (t,x)
\in [0,T] \times \bR^d
 \ \ \ \mbox{if and only if} \ (\ref{cond-H}) \ \mbox{holds}.
 \end{equation}
\end{theorem}

\begin{remark}
{\rm Note that condition (\ref{cond-H}) implies that $d \in
\{1,2,3\}$, since $H<1$. }
\end{remark}

\begin{proof} (a) We will apply Theorem \ref{criterion-cH} to the
function $g_{tx}$. Note that for every fixed $y \in \bR^d$
$$\int_{0}^{T}|g_{tx}(s,y)|^2ds=C\int_{0}^{t}\frac{1}{(t-s)^{d}}
e^{-|y-x|^2/[2(t-s)]}ds=C\int_{1/t}^{\infty}
u^{d-2}e^{-|y-x|^{2}u/2}du< \infty,$$ i.e. $g_{tx}(\cdot,y) \in
L_{2}(0,T)$. Clearly $(\tau,y) \mapsto \cF_{0,T}g_{tx}(\tau,y)$ is
measurable. We now calculate
$$\|g_{tx}\|_{\cH}^{2}:=\alpha_H
c_{H} \int_{\bR}\int_{\bR^d}
|{\cF}_{0,T}g_{tx}(\tau,y)|^{2}|\tau|^{-(2H-1)} dyd\tau.$$

\noindent For this, we write
\begin{eqnarray*}
 \|g_{tx}\|_{\cH}^{2} & =&
\alpha_H c_H \int_{\bR} |\tau|^{-(2H-1)}  \int_{\bR^d}
\left(\int_{0}^{T}e^{-i\tau s}g_{tx}(s,y)ds \right)
 \left(\int_{0}^{T}e^{i\tau r}g_{tx}(r,y)dr \right)
 dy d\tau\\
& = & \alpha_H c_H \int_{\bR} |\tau|^{-(2H-1)} \int_{0}^{t}
\int_{0}^{t} e^{-i\tau (s-r)} I(r,s) dr ds d\tau, \\
& = & \alpha_H \int_{0}^{t} \int_{0}^{t} |s-r|^{2H-2} I(r,s) dr ds
\end{eqnarray*}
where we used Lemma A.1 (Appendix A) for the last equality and we
denoted $$I(r,s):=
\int_{\bR^d}g_{tx}(s,y)g_{tx}(r,y)dy=(2t-s-r)^{-d/2}$$


\noindent We obtain that $$\|g_{tx}\|_{\cH}^{2} = \alpha_H
\int_{0}^{t} \int_{0}^{t}|s-r|^{2H-2}(2t-s-r)^{-d/2}dr ds
 = \alpha_H \int_{0}^{t}
\int_{0}^{t}|u-v|^{2H-2}(u+v)^{-d/2}dr ds.$$

Relation (\ref{nec-suf-cond-H}) follows, since the last integral
is finite if and only if $2H>d/2$. In this case, we have
$\|g_{tx}\|_{|\cH|}=\|g_{tx}\|_{\cH}<\infty$ (since $g_{tx} \geq
0$), and hence $g_{tx} \in |\cH|$.

(b) We will apply Theorem \ref{criterion-cH} to the function
$\varphi=\eta
* \tilde G$, since $\varphi(\cdot,x) \in L_{2}(0,T)$ for every $x \in \bR^d$, by
Lemma \ref{eta*G-in-L2}. By writing
\begin{equation}
\label{Fourier-eta*G} {\cF}_{0,T} \varphi(\tau,x) = \int_{0}^{T}
e^{-i\tau s} \int_{0}^{T-s} \int_{\bR^d} \eta(u+s,y)
G(u,y-x)dyduds,
\end{equation}
we see that $(\tau,x) \mapsto \cF_{0,T}\varphi(\tau,x)$ is
measurable. We now calculate
$$\|\varphi\|_{\cH}^{2}:=\alpha_H c_{H}
\int_{\bR}\int_{\bR^d}|{\cF}_{0,T}\varphi(\tau,x)|^{2}|\tau|^{-(2H-1)}dxd\tau.$$

\noindent Using (\ref{Fourier-eta*G}), we get
\begin{eqnarray}
\nonumber \|\varphi \|_{\cH}^{2}& =&\alpha_H c_H \int_{\bR}
|\tau|^{-(2H-1)} \int_{0}^{T} \int_{0}^{T}
 e^{-i \tau(s-r)}\int_{\bR^d}
\int_{\bR^d} \int_{0}^{T-s} \int_{0}^{T-r}
 \eta(u+s,y) \eta(v+r,z)\\
 \nonumber
& & J(u,v,y,z)dv \, du \,dz\, dy  \,dr \, ds \,d\tau \\
\nonumber &=& \alpha_H \int_{0}^{T} \int_{0}^{T}
|s-r|^{2H-2}\int_{\bR^d} \int_{\bR^d} \int_{0}^{T-s}
\int_{0}^{T-r}
 \eta(u+s,y) \eta(v+r,z)\\
 \label{norm-H-varphi}
& & J(u,v,y,z)dv \, du \,dz\, dy  \,dr \, ds,
\end{eqnarray}

\noindent where we used Lemma \ref{LemmaA1} (Appendix A) for the
second equality and we denoted
$$ J(u,v,y,z):=\int_{\bR^d}
 G(u,y-x) G(v,z-x) dx =\exp\left\{-\frac{|y-z|^{2}}{4(u+v)}\right\}
 (u+v)^{-d/2}.$$

\noindent Clearly
\begin{equation}
\label{bound-for-J} J(u,v,y,z)
\leq (u+v)^{-d/2}
\end{equation}
for every $u \in (0,T-s), v \in (0,T-r)$. Using
(\ref{norm-H-varphi}), (\ref{bound-for-J}) and the fact that $\eta
\in
 \cD((0,T) \times \bR^d)$, we get
\begin{eqnarray*}
\|\varphi \|_{\cH}^{2} &\leq& \alpha_H\int_{0}^{T} \int_{0}^{T}
|s-r|^{2H-2} \int_{\bR^d} \int_{\bR^d} \int_{0}^{T-s}
\int_{0}^{T-r}
|\eta(u+s,y) \eta(v+r,z)| \\
&& (u+v)^{-d/2}  dv \, du \,dz\, dy \, dr  \, ds \\
& \leq &  \alpha_{H}C_{\eta}\int_{0}^{T} \int_{0}^{T} |s-r|^{2H-2}
\int_{0}^{T-s} \int_{0}^{T-r}(u+v)^{-d/2} dv \, du \, dr \,ds \\
 & \leq & \alpha_{H}C_{\eta}\int_{0}^{T} \int_{0}^{T}(u+v)^{-d/2}
 (T-u)^{2H} dv \, du,
\end{eqnarray*}
where for the last inequality we used Fubini's theorem and the
fact that $\int_{0}^{T-u}\int_{0}^{T-v}|s-r|^{2H-2} dr \,
ds=R_{H}(T-u,T-v)=[(T-u)^{2H}+(T-v)^{2H}-(u-v)^{2H}]/2$. The last
integral is clearly finite if $2H>d/2$, i.e. $H>d/4$. \qed
\end{proof}

\vspace{3mm}

Under the conditions of Theorem \ref{H-psi*Gtilde}, $F(g_{tx})$
and $F(\eta
* \tilde G)$ are well-defined for every $(t,x)$, respectively for
every $\eta \in {\cD}((0,T) \times \mathbb{R}^{d})$ and we can
introduce the following definition:

\begin{definition}
{\rm a) The process $\{v(t,x);t \in [0,T], x \in \bR^{d}\}$
defined by
\begin{equation}
\label{solution-process-H} v(t,x):=F(g_{tx})=\int_{0}^{T}
\int_{\bR^{d}} G(t-s,x-y)F(ds,dy)
\end{equation}
is called the {\bf process solution} of the stochastic heat
equation (\ref{heat-eq-H}), with vanishing initial conditions.

b) The process $\{v(\eta); \eta \in {\cD}((0,T) \times
\mathbb{R}^{d})\}$ defined by
$$v(\eta):  = F(\eta * \tilde G)=\int_{0}^{T} \int_{\bR^{d}} \left( \int_{\bR_{+}}
\int_{\bR^{d}}\eta(t+s,x+y)G(s,y)dyds \right) F(dt,dx)$$ is called
the {\bf distribution-valued solution} of the stochastic heat
equation (\ref{heat-eq-H}), with vanishing initial conditions.}
\end{definition}

\begin{lemma}
\label{L2-cont-v} The process $\{v(t,x);t \in [0,T], x \in
\bR^{d}\}$ is $L^{2}(\Omega)$-continuous.
\end{lemma}

\begin{proof}
We first prove the continuity in $t$. We have
\begin{eqnarray*}
\lefteqn{E|v(t+h,x)-v(t,x)|^{2}= E\left|\int_{0}^{t+h}\int_{\bR^d}
g_{t+h,x}(s,y)F(ds,dy)-
\int_{0}^{t}\int_{\bR^d}g_{tx}(s,y)F(ds,dy)
\right|^{2} } \\
&=&E\left|\int_{0}^{t}\int_{\bR^d}
(g_{t+h,x}(s,y)-g_{tx}(s,y))F(ds,dy)+
\int_{t}^{t+h}\int_{\bR^d}g_{t+h,x}(s,y)F(ds,dy) \right|^{2} \\
& \leq & 2E\left|\int_{0}^{t} \int_{\bR^d}
(g_{t+h,x}(s,y)-g_{tx}(s,y))F(ds,dy)\right|^{2}+
2E\left|\int_{t}^{t+h} \int_{\bR^d}
g_{t+h,x}(s,y)F(ds,dy)\right|^{2} \\
&:=& 2I_{1}(h)+2I_{2}(h).
\end{eqnarray*}

\noindent We first treat the term $I_{1}(h)$. Note that
\begin{eqnarray*}
I_{1}(h)&=&\|F((g_{t+h,x}-g_{tx})1_{[0,t]})\|^2_{L_{2}(\Omega)}
= \| (g_{t+h,x}-g_{tx})1_{[0,t]}\|_{\cH}^2= \\
&=& \alpha_H \int_{0}^{t}\int_{0}^{t}
\int_{\bR^d}(g_{t+h,x}-g_{tx})(u,y)|u-v|^{2H-2}
(g_{t+h,x}-g_{tx})(v,y)dy \, dv \, du.
\end{eqnarray*}

\noindent The continuity of $G(t,x)$ with respect to $t$ shows
that the integrand converges to zero as $h\to 0$. By applying the
dominated convergence theorem (using the fact that
$\|g_{tx}\|_{\cH}<\infty$), we conclude that $I_{1}(h) \to 0$.

For the term $I_2(h)$, we have
\begin{eqnarray*}
I_{2}(h)&=&\|F(g_{t+h,x}1_{[t,t+h]})\|^2_{L_{2}(\Omega)}=\|
g_{t+h,x}1_{[t,t+h]}\|_{\cH}^2 =\\
 &=& \alpha _{H}\int _{t}^{t+h}
\int _{t}^{t+h} \int_{\mathbb{R}^{d}}g_{t+h, x}(u,y) g_{t+h,
x}(v,y)|u-v|^{2H-2} dy \, dv \, du \\
&=&\alpha _{H}\int _{t-h}^{t} \int _{t-h}^{t}\int_{\bR^{d}}
g_{tx}(u',y) g_{tx}(v',y)|u'-v'|^{2H-2} dy \, dv' \, du.
\end{eqnarray*}

\noindent Since $1_{(t-h,t)}1_{(t-h,t)} \to 0$ as $h\to 0$ and
$\|g_{tx}\|_{\cH} <\infty$, we conclude that $I_2(h) \to 0$, by
the dominated convergence theorem.

We now prove the continuity in $x$. We have
\begin{eqnarray*}
E|v(t,x+h)-v(t,x)|^{2}&=& E\left|\int_{0}^{t}\int_{\bR^d}
(g_{t,x+h}(s,y)-g_{tx}(s,y))F(ds,dy)\right|^{2} =
\|(g_{t,x+h}-g_{tx})1_{[0,t]}\|_{\cH}^2 \\
&=&\alpha_H \int_{0}^{t}
\int_{0}^{t}(g_{t,x+h}-g_{tx})(u,y)||u-v|^{2H-2}
(g_{t,x+h}-g_{tx})(v,y) dy \, dv \, du.
\end{eqnarray*}

\noindent By the continuity in $x$ of the function  $G(t,x)$, the
integrand converges to 0 as $h\to 0$. By the dominated convergence
theorem, we conclude that $E|v(t,x+h)-v(t,x)|^{2} \to 0$ as $h \to
0$. \qed
\end{proof}

\vspace{3mm}

\begin{theorem}
\label{existence-theorem-H} Let $\{v(\eta); \eta \in {\cD}((0,T)
\times \mathbb{R}^{d})\}$ be the distribution-valued solution of
the stochastic heat equation (\ref{heat-eq-H}).

In order that there exists a jointly measurable and locally
mean-square bounded process $Y=\{Y(t,x);t \in [0,T],x \in
\bR^{d}\}$ such that
\begin{equation}
\label{v(eta)=integral(v)} v(\eta)=\int_{0}^{T} \int_{\bR^{d}}
Y(t,x)\eta(t,x) dxdt \ \ \ \forall \eta \in {\cD}((0,T) \times
\bR^{d}) \ \ \ a.s.
\end{equation}
 it is necessary and
sufficient that (\ref{cond-H}) holds.  In this case, $Y$ is a
modification of the process $v=\{v(t,x);t \in [0,T],x \in
\bR^{d}\}$ defined by (\ref{solution-process-H}).
\end{theorem}

\begin{proof} The necessity part is similar to the proof of Theorem 11,
\cite{Da1}. Suppose that there exists a jointly measurable and
locally mean-square bounded process $Y=\{Y(t,x);t \in [0,T],x \in
\bR^{d}\}$ such that (\ref{v(eta)=integral(v)}) holds.

Let $(t_0,x_0) \in [0,T] \times \bR^d$ be fixed. Let
$\lambda,\psi$ be nonnegative functions such that $\lambda \in
\cD(0,T), \psi \in \cD(\bR^d)$ and $\int_0^T
\lambda(t)dt=\int_{\bR^d}\psi(x)dx=1$. Set
$\lambda_n(t)=n\lambda(nt)$ and $\psi_n(x)=n^d\psi (nx)$. Define
$\eta_n(t,x)=\lambda_n(t-t_0) \psi_n(x-x_0)$.

We calculate $E|v(\eta_n)|^{2}$ in two ways. First, using
(\ref{v(eta)=integral(v)}), we get
$$E|v(\eta_n)|^2=\int_{[0,T] \times \bR^d}
E|Y(t,x)Y(s,y)| \ \eta_n(t,x)\eta_{n}(s,y)dydxdsdt.$$

\noindent Using Lebesgue differentiation theorem (see Exercise 2,
Chapter 7, \cite{wheeden-zygmund77}), we get
\begin{equation}
\label{expect-Ev(eta-n)-1}
 \lim_n
E|v(\eta_n)|^{2}=E|Y(t_0,x_0)|^{2}.
\end{equation}

\noindent Secondly,
\begin{equation}
\label{expect-Ev(eta-n)-2} E|v(\eta_n)|^{2}=E|F(\eta_n * \tilde
G)|^2=\|\eta_n * \tilde G \|_{\cH}^{2} =\int_{\bR^d}
\int_{\bR}|\cF_{0,T}(\eta_n * \tilde
G)(\tau,x)|^{2}|\tau|^{-(2H-1)}d \tau dx.\nonumber
\end{equation}

\noindent We claim that, for every $\tau \in \bR, x \in \bR^d$, we
have: (see Appendix B for the proof)
\begin{equation}
\label{Theorem11-Dalang-claim} \lim_n \cF_{0,T} (\eta_n * \tilde
G)(\tau,x)=\cF_{0,T} g_{t_0 x_0}(\tau,x).
\end{equation}

\noindent Using Fatou's lemma, (\ref{expect-Ev(eta-n)-1}) and
(\ref{expect-Ev(eta-n)-2}), we get
\begin{eqnarray*} \|g_{t_0 x_0} \|_{\cH}^{2}&=& \int_{\bR^d}
\int_{\bR}|\cF_{0,T} g_{t_0
x_0}(\tau,x)|^{2}|\tau|^{-(2H-1)}d \tau dx \\
&\leq &\liminf_{n} \int_{\bR^d} \int_{\bR} |\cF_{0,T} (\eta_n *
\tilde G) (\tau,x)|^{2}|\tau|^{-(2H-1)}d \tau dx\\
 &=&\liminf_n
E|v(\eta_n)|^{2}=E|Y(t_0,x_0)|^{2}<\infty, \end{eqnarray*} which
forces $H>d/4$, by virtue of (\ref{nec-suf-cond-H}).

We now prove the sufficiency part. By Lemma \ref{L2-cont-v}, we
know that the process $v$ defined by (\ref{solution-process-H}) is
continuous in $L^{2}(\Omega)$. Hence, it is continuous in
probability and by Theorem IV.30, \cite{dellacherie-meyer75},
there exists a jointly measurable modification $Y$ of $v$. We have
\begin{equation}
\label{cov-Y}
E(Y(t,x)Y(s,y))=E(v(t,x)v(s,y))=E(F(g_{tx})F(g_{sy}))=\langle
g_{tx},g_{sy} \rangle_{\cH}
\end{equation}
and
\begin{equation}
\label{cov-v-Y} E(v(\eta)Y(t,x))=E(v(\eta)v(t,x))=E(F(\eta *
\tilde G)F(g_{tx}))=\langle \eta * \tilde G, g_{tx} \rangle_{\cH}.
\end{equation}

To prove that (\ref{v(eta)=integral(v)}) holds, we will show that
$$E\left|\int_{0}^{T} \int_{\bR^{d}}
Y(t,x)\eta(t,x) dxdt \right|^{2}=E\left(v(\eta) \int_{0}^{T}
\int_{\bR^{d}} Y(t,x)\eta(t,x) dxdt \right)=E|v(\eta)|^{2}.$$

\noindent By Fubini's theorem and (\ref{cov-Y}), we get
\begin{eqnarray*}
\lefteqn{E\left|\int_{0}^{T} \int_{\bR^{d}} Y(t,x)\eta(t,x) dxdt
\right|^{2} =\int_{([0,T] \times \bR^d)^{2}} \eta(t,x)\eta(s,y)
\langle
g_{tx},g_{sy} \rangle_{\cH} dy \, dx \, ds \, dt =} \\
&& \int_{([0,T] \times \bR^d)^{2}} \eta(t,x)\eta(s,y) \left(
\int_{\bR} \int_{\bR^d} \cF_{0,T}g_{tx}(\tau,z)
\overline{\cF_{0,T}
g_{sy}(\tau,z)} \ |\tau|^{-(2H-1)} dz d \tau \right) dydx ds dt =\\
&& \int_{\bR} \int_{\bR^d} |\cF_{0,T} (\eta * \tilde
G)(\tau,z)|^{2}|\tau|^{-(2H-1)} dz \, d \tau=\|\eta * \tilde G
\|_{\cH}^{2}=E|F(\eta * \tilde G)|^{2}=E|v(\eta)|^{2}.
\end{eqnarray*}

\noindent On the other hand, using Fubini's theorem and
(\ref{cov-v-Y}), we get
\begin{eqnarray*}
\lefteqn{E\left(v(\eta) \int_{0}^{T} \int_{\bR^{d}}
Y(t,x)\eta(t,x) dxdt \right)=\int_{[0,T] \times \bR^d}  \eta(t,x)
\langle \eta *
\tilde G, g_{tx}\rangle_{\cH} dxdt =} \\
&& \int_{[0,T] \times \bR^d} \eta(t,x)
\left(\int_{\bR}\int_{\bR^d} \cF_{0,T}(\eta * \tilde G)(\tau,z)
\overline{\cF_{0,T}
g_{tx}(\tau,z)} \ |\tau|^{-(2H-1)}  dz d \tau\right) dx dt =\\
&& \int_{\bR}  \int_{\bR^d} |\cF_{0,T} (\eta * \tilde
G)(\tau,z)|^{2}|\tau|^{-(2H-1)} dz d \tau=E|v(\eta)|^{2}.
\end{eqnarray*}

\qed
\end{proof}

\section{The Fractional-Colored Noise}
\label{sect-fract-color}

In this section we examine the process-solution and the
distribution-solution of the stochastic heat equation driven by a
Gaussian noise which is {\em fractional in time} and {\em colored
in space}. Most of the results of this section are obtained by
mixing some colored spatial techniques with the fractional
temporal techniques of Section \ref{sect-fract-white}. The results
of this section can therefore be viewed as generalizations of the
results of Section \ref{sect-fract-white}. The details are highly
non-trivial.

The structure of this section is similar to that of Section
\ref{sect-fract-white}. We first describe the spaces of
deterministic integrands, then we introduce the Gaussian noise and
the associated stochastic integral, and finally we examine the
solution of the stochastic heat equation driven by this type of
noise.

\subsection{Spaces of Deterministic Integrands}

We begin by introducing the space of deterministic integrands on
$\bR^d$.

We say that a function $f:\bR \to \bR$ is a {\em spatial
covariance function}, if it is the Fourier transform of a tempered
measure $\mu$ on $\bR^d$, i.e. $f(x)=\int_{\bR^d}e^{-i \xi \cdot
x}\mu(d \xi)$.

Let $\cP(\bR^d)$ be the completion of $\cD(\bR^d)$ with respect to
the inner product
$$\langle \varphi,\psi
\rangle_{\cP(\bR^d)}=\int_{\bR^{d}} \int_{\bR^{d}}
\varphi(x)f(x-y) \psi(y)dy
dx=\int_{\bR^d}\cF_{2}\varphi(\xi)\overline{\cF_{2}\psi(\xi)}\mu(d
\xi),$$ where $\cF_{2}\varphi(\xi):=\int_{\bR^d}e^{-i \xi
x}\varphi(x)dx$ denotes the Fourier transform with respect to the
$x$-variable.

Equivalently, we can say that $\cP(\bR^d)$ is the completion of
$\cE(\bR^d)$ with respect to $\langle \cdot, \cdot
\rangle_{\cP(\bR^d)}$, where $\cE(\bR^d)$ is the space of all
linear combinations of indicator functions $1_{A}(x), A \in
\cB_{b}(\bR^d)$)

The basic example of spatial covariance function is
$f(x)=\delta(x)$, which gives rise to the spatial white noise.
More interesting covariance structures are provided by potential
analysis. Here are some examples (see e.g. p.149-151,
\cite{folland95}, or p.117-132, \cite{stein70}; our constants are
slightly different than those given in these references, since our
definition of the Fourier transform does not have the $2\pi$
factor):

\begin{example}
\label{riesz} {\rm  The Riesz kernel of order $\alpha$:
$$f(x)=R_{\alpha}(x):=
\gamma_{\alpha,d}|x|^{-d+\alpha}, \quad 0<\alpha<d,$$
 where $\gamma_{\alpha,d}=2^{d-\alpha} \pi^{d/2}
 \Gamma((d-\alpha)/2)/\Gamma(\alpha/2)$. In this case, $\mu(d\xi)=
 |\xi|^{-\alpha}d \xi$. }
\end{example}

\begin{example}
\label{bessel}{\rm The Bessel kernel of order $\alpha$:
$$f(x)=B_{\alpha}(x):=\gamma'_{\alpha}\int_{0}^{\infty}w^{(\alpha-d)/2-1}
e^{-w}e^{-|x|^{2}/(4w)} dw, \quad \alpha>0,$$ where
$\gamma'_{\alpha}=(4\pi)^{\alpha/2}\Gamma(\alpha/2)$. In this
case, $\mu(d\xi)=(1+|\xi|^{2})^{-\alpha/2} d \xi$ and $\cP(\bR^d)$
coincides with $\cH^{-\alpha/2}(\bR^d)$, the fractional Sobolev
space of order $-\alpha/2$; see e.g. p.191, \cite{folland95}. }
\end{example}

\begin{example} {\rm The heat kernel
$$f(x)=G_{\alpha}(x):= \gamma_{\alpha,d}''
e^{-|x|^{2}/(4\alpha)}, \quad \alpha>0,$$ where
$\gamma_{\alpha,d}''=(4 \pi \alpha)^{-d/2}$. In this case,
$\mu(d\xi)=e^{-\pi^2 \alpha |\xi|^2} d \xi$.}
\end{example}

\begin{example} {\rm The Poisson kernel
$$f(x)=P_{\alpha}(x):=\gamma_{\alpha,d}'''
(|x|^{2}+\alpha^{2})^{-(d+1)/2}, \quad \alpha>0,$$ where
$\gamma_{\alpha,d}'''=\pi^{-(d+1)/2}\Gamma((d+1)/2)\alpha$. In
this case, $\mu(d\xi)=e^{-4\pi^2 \alpha |\xi|} d \xi$. }
\end{example}

\begin{remark}
{\rm  The space $\cP$ defined as the completion of $\cD((0,T)
\times \bR^d)$ (or the completion of $\cE$) with respect to the
inner product
$$\langle \varphi,\psi \rangle_{\cP}=\int_{0}^{T} \int_{\bR^{d}}
\int_{\bR^{d}} \varphi(t,x)f(x-y) \psi(t,y)dy dx dt=\int_{0}^{T}
\langle \varphi(t, \cdot), \psi(t, \cdot)
\rangle_{\cP(\bR^d)}dt.$$ has been studied by several authors in
connection with a Gausssian noise which is white in time and
colored in space. One can prove that $\cP \subset L_{2}((0,T);
\cP(\bR^d))$; see e.g. \cite{Da1}, or \cite{Balan}. }
\end{remark}

In what follows, we need to extend the definition of $\cP(\bR^d)$
to allow for complex-valued functions. More precisely, let
$\cD_{\bC}(\bR^d)$ be the space of all infinitely differentiable
functions $\varphi: \bR^d \to \bC$ with compact support, and
$\cP_{\bC}(\bR^d)$ be the completion of $\cD_{\bC}(\bR^d)$ with
respect to
$$\langle \varphi, \psi \rangle_{\cP_{\bC}
(\bR^d)}=\int_{\bR^{d}} \int_{\bR^{d}} \varphi(x)f(x-y) \psi(y)dy
dx.$$

\noindent Since $\cD(\bR^d) \subset \cD_{\bC}(\bR^d)$ and $\langle
\varphi, \psi \rangle_{\cP(\bR^d)}=\langle \varphi, \psi
\rangle_{\cP_{\bC}(\bR^d)}$ for every $\varphi,\psi \in
\cD(\bR^d)$, we conclude that $$\cP(\bR^d) \subset
\cP_{\bC}(\bR^d).$$

We are now introducing the space of deterministic integrands
associated with a Gaussian noise which is fractional in time and
colored in space. This space seems to be new in the literature.
More precisely, let $\cH \cP$ be the completion of $\cD((0,T)
\times \bR^d)$ with respect to the inner product
\begin{eqnarray*}
\langle \varphi, \psi \rangle_{\cH \cP}&=& \alpha_{H}\int_{0}^{T}
\int_{0}^{T} \int_{\bR^d}\int_{\bR^d}
\varphi (u,x) |u-v|^{2H-2}f(x-y) \psi (u,y)dy \, dx \, dv \, du \\
&=&\alpha_{H}c_{H}  \int_{\bR}|\tau|^{-(2H-1)}
\int_{\bR^d}\int_{\bR^d} f(x-y) \cF_{0,T}\varphi(\tau,x)
\overline{\cF_{0,T}\psi(\tau,y)}
dy \, dx \, d \tau \\
& = & \alpha_{H}c_{H}  \int_{\bR} |\tau|^{-(2H-1)} \langle
\cF_{0,T}\varphi(\tau,\cdot), \overline{\cF_{0,T}\psi(\tau,\cdot)}
\rangle_{\cP_{\bC}(\bR^d)}  d \tau,
\end{eqnarray*}
where the second equality follows by Lemma \ref{LemmaA1} (Appendix
A). In particular,
$$\|\varphi \|_{\cH \cP}^2=\alpha_{H}c_{H}  \int_{\bR}  \|
\cF_{0,T}\varphi(\tau,\cdot)\|_{\cP_{\bC}(\bR^d)}^2
|\tau|^{-(2H-1)} d \tau,$$ or equivalently $\|\varphi \|_{\cH
\cP}^2=\int_{\bR}
\|\cF_{0,T}\varphi(\tau,\cdot)\|_{\cP_{\bC}(\bR^d)}^2
\lambda_{H}(d \tau)$, where $\lambda_H(d \tau)=\alpha_H
c_{H}|\tau|^{-(2H-1)}d \tau$.

One can prove that $\cH \cP$ is also the completion of $\cE$ with
respect to the inner product
$$\langle 1_{[0,t] \times A}, 1_{[0,s] \times B} \rangle_{\cH
\cP}=R_{H}(t,s) \langle 1_{A},1_{B} \rangle_{\cP(\bR^d)}.$$

\noindent Clearly $|\cH \cP| \subset \cH  \cP$, where $|\cH
\cP|=\{\varphi:[0,T]\times \bR^d \ \mbox{measurable}; \  \|\varphi
\|_{|\cH \cP|}<\infty\}$ and  $$\|\varphi \|_{|\cH
\cP|}^2:=\int_{0}^{T} \int_{0}^{T} \int_{\bR^d}
\int_{\bR^d}|\varphi(u,x)| |\varphi(v,y)| |u-v|^{2H-2}f(x-y)  dy
\, dx \, dv \, du.$$

\begin{remark}
{\rm Using Fubini's theorem, we have the following alternative
expression for calculating $\langle \varphi, \psi \rangle_{\cH
\cP}$: for every $\varphi,\psi \in \cD((0,T) \times \bR^d)$, we
have
\begin{eqnarray*}
\langle \varphi, \psi \rangle_{\cH \cP}&=& \alpha_{H}\int_{0}^{T}
\int_{0}^{T} |u-v|^{2H-2} \int_{\bR^d}
\cF_2\varphi (u,\xi) \overline{\cF_2\psi(u,\xi)} \mu(d \xi)dv du \\
&=&\alpha_{H} \int_{\bR^d}\int_{0}^{T} \int_{0}^{T}
\cF_{2}\varphi(u,\xi) |u-v|^{2H-2}\overline{\cF_2\psi(v,\xi)} dv
du
\mu(d \xi) \\
& = & \int_{\bR^d} \langle \cF_{2}\varphi(\cdot,\xi),
\overline{\cF_2\psi(\cdot,\xi)} \rangle_{\cH_{\bC}(0,T)}  \mu(d
\xi),
\end{eqnarray*}
where $\cH_{\bC}(0,T)$ denotes the completion of $\cD_{\bC}(0,T)$
with respect to the inner-product $\langle \cdot, \cdot
\rangle_{\cH_{\bC}(0,T)}$ defined similarly to $\langle \cdot,
\cdot \rangle_{\cH(0,T)}$. In particular, $\|\varphi \|_{\cH
\cP}^2=\int_{\bR^d} \| \cF_{2}\varphi(\cdot,\xi)
\|_{\cH_{\bC}(0,T)}^2 \mu(d \xi)$.
 This expression will not be used in the present paper.}
\end{remark}

In this new context, the next theorem gives us a useful criterion
for verifying that a function $\varphi$ lies in $\cH \cP$. To
prove this theorem, we need the following lemma, generalizing
Lemma 5.1, \cite{PT1}.

\begin{lemma}
\label{indic-approx-Fourier2} For every $A \in \cB_b(\bR^d)$,
there exists a sequence $(g_n)_{n} \subset \cE(\bR^d)$ such that
$$\int_{\bR^d}|1_{A}(\xi)-\cF_{2}g_{n}(\xi)|^{2}\mu(d\xi) \to 0.$$
\end{lemma}

\begin{proof}
Let $(\phi_n)_n \subset \cD(\bR^d)$ be such that $\phi_n(\xi) \to
1_{A}(\xi)$ uniformly and ${\rm supp} \ \phi_n \subset K,\ \forall
n$, where $K \subset \bR^d$ is a compact (we may take
$\phi_n=1_{A}*\eta_n$, where $\eta_n(x)=n^d\eta(nx)$ and $\eta \in
\cD(\bR^d)$, with $\int_{\bR^d} |\eta(x)|dx=1$). Using the
dominated convergence theorem and the fact that $\mu$ is locally
finite, we get $\int_{\bR^d}|\phi_n(\xi)-1_{A}(\xi)|^2\mu(d\xi)
\to 0$.

Let $\psi_n \in \cS(\bR^d)$ be such that $\cF_2 \psi_n=\phi_n$.
Then
\begin{equation}
\label{Fourier-psi-n-1A}
\int_{\bR^d}|\cF_{2}\psi_n(\xi)-1_{A}(\xi)|^2\mu(d\xi) \to 0.
\end{equation}

\noindent Note that $\int_{\bR^d}|\cF_{2}\psi_n(\xi)|^2
\mu(d\xi)=\int_{\bR^d}|\phi_n(\xi)|^2 \mu(d \xi)<\infty$. By Lemma
\ref{LemmaC1} (Appendix C), it follows that $\psi_n \in
\cP(\bR^d)$. Since $\cP(\bR^d)$ is the completion of $\cE(\bR^d)$
with respect to $\| \cdot \|_{\cP(\bR^d)}$, there exists a
sequence $(g_{n})_n \subset \cE(\bR^d)$ such that
\begin{equation}
\label{Fourier-psi-n-g-n} \|\psi_n-g_{n}\|_{\cP(\bR^d)}^{2}=
\int_{\bR^d}|\cF_{2}\psi_n(\xi)-\cF_2g_{n}(\xi)|^2\mu(d\xi) \to 0.
\end{equation}

\noindent The conclusion follows from (\ref{Fourier-psi-n-1A}) and
(\ref{Fourier-psi-n-g-n}). \qed
\end{proof}

\vspace{3mm}

\begin{theorem}
\label{criterion-cHcP} Let $\varphi :[0,T] \times \bR^d
\rightarrow \bR$ be a function which satisfies the following
conditions:

(i) $\varphi(\cdot,x) \in L_{2}(0,T)$ for every $x \in \bR^d$;

(ii) $(\tau,x) \mapsto \cF_{0,T} \varphi (\tau,\cdot)$ is
measurable;

(iii) $\int_{\bR} \int_{\bR^d} \int_{\bR^d} \cF_{0,T}
\varphi(\tau,x) f(x-y) \overline{\cF_{0,T}
\varphi(\tau,y)}|\tau|^{-(2H-1)} dydx d\tau<\infty$.

\noindent Then $\varphi \in \cH \cP$.
\end{theorem}

\begin{proof}
The proof follows the same lines as the proof of Theorem
\ref{criterion-cH}. The details are quite different though. Let
\begin{eqnarray*}
\tilde \Lambda& =& \{\varphi:[0,T] \times \bR^d \to \bR; \
\varphi(\cdot,x) \in L_{2}(0,T) \ \forall x, (\tau,x) \mapsto
\cF_{0,T} \varphi (\tau,x) \ \mbox{is measurable, and} \\
& & \|\varphi \|_{\tilde \Lambda}^{2}:=c_{1} \int_{\bR}
\int_{\bR^d} \int_{\bR^d}\cF_{0,T} \varphi(\tau,x)f(x-y)
\overline{\cF_{0,T} \varphi(\tau,y)} |\tau|^{-(2H-1)} dy dx
d\tau<\infty\}
\\
\Lambda&=&\{\varphi:[0,T] \times \bR^d \to \bR; \ \|\varphi
\|_{\Lambda}^{2}:=c_{2}\int_{0}^{T} \int_{\bR^d} \int_{\bR^d}
I_{T-}^{H-1/2}(u^{H-1/2}\varphi(u,x))(s) \cdot f(x-y) \\
& & \times I_{T-}^{H-1/2}(u^{H-1/2}\varphi(u,y))(s) \cdot
s^{-(2H-1)} dy dx ds <\infty\}
\end{eqnarray*}
where $c_{1}=\alpha_Hc_{H}$ and
$c_2=\{c_{H}^*\Gamma(H-1/2)\}^{2}$. The fact that
\begin{equation}
\label{cP-tildeLambda-in-Lambda} \tilde \Lambda \subset \Lambda
\quad \mbox{and} \quad \|\varphi \|_{\tilde \Lambda}=\|\varphi
\|_{\Lambda}, \quad \forall \varphi \in \tilde \Lambda.
\end{equation}
follows as in the proof of Theorem \ref{criterion-cH}: let
$\varphi \in \tilde \Lambda$ is arbitrary; using the fact
$K_{H}^{*}$ is an isometry from $\cH(0,T)$ to $L_{2}(0,T)$, we get
$\langle \varphi(\cdot,x), \varphi(\cdot,y)
\rangle_{\cH(0,T)}^{2}=\langle K_{H}^{*}\varphi (\cdot,x),
K_{H}^{*} \varphi(\cdot,y) \rangle_{L_{2}(0,T)}^{2}$ for all $x,y
\in \bR^d$. Multiplying by $f(x-y)$, integrating with respect to
$dx \, dy$ and using Fubini's theorem, we get $\|\varphi
\|_{\tilde \Lambda}=\|\varphi \|_{\Lambda}<\infty$.

Next we prove that
\begin{equation}
\label{cP-E-dense-in-Lambda} \cE \ \mbox{is dense in} \ \Lambda \
\mbox{with respect to} \ \|\cdot \|_{\Lambda}.
\end{equation}

\noindent The proof of the theorem will follow from
(\ref{cP-tildeLambda-in-Lambda}) and (\ref{cP-E-dense-in-Lambda}),
as in the case of Theorem \ref{criterion-cH}.

To prove (\ref{cP-E-dense-in-Lambda}), let $\varphi \in \Lambda$
and $\varepsilon>0$ be arbitrary. Let
$\lambda_{H}(s)=s^{-(2H-1)}ds$ and \linebreak
$a(s,x)=I_{T-}^{H-1/2}(u^{H-1/2}\varphi(u,x))(s)$. First, we claim
that there exists $g \in \cE$ such that
\begin{equation}
\label{cHcP-first-approx-IT-varphi} I_1:=\int_{0}^{T}\int_{\bR^d}
\int_{\bR^d}[a(s,x)-g(s,x)]f(x-y)[a(s,y)-g(s,y)]dy dx
\lambda_{H}(ds) <\varepsilon.
\end{equation}

\noindent To see this, note that
$$\|\varphi\|_{\Lambda}^{2}=c_{2} \int_{0}^{T}\int_{\bR^d}
\int_{\bR^d}a(s,x)f(x-y)a(s,y)dy dx \lambda_{H}(ds)
=c_{2}\int_{0}^{T}\int_{\bR^d}
|\cF_{2}a(s,\xi)|\mu(d\xi)\lambda_{H}(ds)<\infty,$$ i.e. the map
$(s,\xi) \mapsto \cF_{2}a(s,\xi)$ belongs to $L_{2}((0,T) \times
\bR^d,d\lambda_H \times d\mu)$. Hence, there exists a simple
function
$h(s,\xi)=\sum_{j=1}^{m}\beta_{j}1_{[\gamma_j,\delta_j)}(s)1_{B_j}(\xi)$
on $(0,T) \times \bR^d$, with $\beta_{j} \in \bR$,
$0<\gamma_j<\delta_j<T$ and $B_{j} \subset \bR^d$ Borel sets, such
that
\begin{equation}
\label{approx-Fourier2}
\int_{0}^{T}\int_{\bR^d}|\cF_{2}a(s,\xi)-h(s,\xi)|^{2}\mu(d\xi)\lambda_{H}(ds)
< \varepsilon/4.
\end{equation}

\noindent By Lemma \ref{indic-approx-Fourier2}, for every $j=1,
\ldots, m$, there exists $g_{j} \in \cE(\bR^d)$ such that
\begin{equation}
\label{indic-approx-Fourier2-level-j}
\int_{\bR^d}|1_{B_j}(\xi)-\cF_{2}g_{j}(\xi)|\mu(d\xi)<\varepsilon/(4D_{h}),
\end{equation}
where we chose
$D_{h}=m\sum_{j=1}^{m}\beta_{j}^{2}\lambda_{H}([\gamma_j,\delta_j))$.
Define $g(s,x)=\sum_{j=1}^{m}\beta_{j}
1_{[\gamma_j,\delta_j)}(s)g_{j}(x)$. Clearly $g \in \cE$. Using
(\ref{approx-Fourier2}) and (\ref{indic-approx-Fourier2-level-j}),
we get
$I_{1}=\int_{0}^{T}\int_{\bR^d}|\cF_{2}a(s,\xi)-\cF_{2}g(s,\xi)|^{2}\mu(d\xi)
\lambda_{H}(ds)<\varepsilon$, which concludes the proof of
(\ref{cHcP-first-approx-IT-varphi}).

We claim now that there exists a function $l \in \cE$ such that
\begin{equation}
\label{cHcP-second-approx-IT-varphi} I_2:=\int_{0}^{T}
\int_{\bR^d} \int_{\bR^d} [g(s,x)-b(s,x)]f(x-y)[g(s,y)-b(s,y)]dy
dx \lambda_{H}(ds)<\varepsilon,
\end{equation}
where $b(s,x)=I_{T-}^{H-1/2}(u^{H-1/2}l(u,x))(s)$. To see this,
suppose that $g=\sum_{k=1}^{n}b_k 1_{[c_k,d_k)}(s)1_{A_k}(x)$ for
some $0<c_k<d_k<T$ and $A_k \subset \bR^d$ Borel sets. By relation
(8.1) of \cite{PT}, there exists an elementary function $l_{k} \in
\cE(0,T)$ such that
\begin{equation}
\label{rel8-1-PT} \int_{0}^{T}
[1_{[c_k,d_k)}(s)-I_{T-}^{H-1/2}(u^{H-1/2}l_k(u))(s)]^{2}\lambda_{H}(ds)<
\varepsilon/C_g,
\end{equation}
 where we chose $C_g:=\|
\sum_{k=1}^{n}b_{k}1_{A_k} \|_{\cP(\bR^d)}^{2}$. Let
$l(s,x)=\sum_{k=1}^{n}b_k l_{k}(t)1_{A_k}(x) \in \cE$ and note
that
\begin{eqnarray*}
I_2 &=&\sum_{k,j=1}^{n}b_{k} b_{j}\int_{0}^{T} \int_{\bR^d}
\int_{\bR^d}   [1_{[c_k,d_k)}(s)-
I_{T-}^{H-1/2}(u^{H-1/2}l_k(u))(s)]1_{A_k}(x) f(x-y) \\
& & [1_{[c_j,d_j)}(s)-I_{T-}^{H-1/2}(u^{H-1/2}l_j(u))(s)]
1_{A_j}(y)dy dx \lambda_{H}(ds) \\
& \leq & \sum_{k,j=1}^{n}b_k b_j \langle 1_{A_k},1_{A_j}
\rangle_{\cP(\bR^d)} \cdot 2
\left(\frac{\varepsilon}{C_g}+\frac{\varepsilon}{C_g}\right)
=\varepsilon
\end{eqnarray*}
(we used (\ref{rel8-1-PT}) and the fact that $ab \leq
2(a^2+b^2)$). The proof of (\ref{cHcP-second-approx-IT-varphi}) is
complete.

Finally, we claim that from (\ref{cHcP-first-approx-IT-varphi})
and (\ref{cHcP-second-approx-IT-varphi}), we get
$$I:=\int_{0}^{T}\int_{\bR^d}[a(s,x)-b(s,x)]f(x-y)[a(s,y)-b(s,y)]
dydx \lambda_{H}(ds) <4\varepsilon,$$ i.e. $\|\varphi-l
\|_{\Lambda}^{2}<4 \varepsilon c_2$, which will conclude the proof
of (\ref{cP-E-dense-in-Lambda}). To see this, note that
$I=\sum_{k=1}^{4}I_{k}$, where
\begin{eqnarray*}
I_{3}&:=&-\int_{0}^{T} \int_{\bR^d} \int_{\bR^d}
[a(s,x)-g(s,x)]f(x-y)[b(s,y)-g(s,y)]dy
dx\lambda_{H}(ds) \\
I_{4}&:=&-\int_{0}^{T} \int_{\bR^d} \int_{\bR^d}
[a(s,y)-g(s,y)]f(x-y)[b(s,x)-g(s,x)]dy dx\lambda_{H}(ds).
\end{eqnarray*}
The fact that $|I_3|<\varepsilon$ and $|I_4|<\varepsilon$ follows
from the Cauchy-Schwartz inequality in $\cP(\bR^d)$,
(\ref{cHcP-first-approx-IT-varphi}) and
(\ref{cHcP-second-approx-IT-varphi}). \qed
\end{proof}

\vspace{3mm}

It is again possible to describe the space $\cH \cP$ using the
transfer operator.  Define the transfer operator on $\cE$ by the
same formula (\ref{def-K-H*}). Note that in this case we have
\begin{eqnarray*}
\langle K_{H}^* 1_{[0,t] \times A}, K_{H}^*1_{[0,s] \times B}
\rangle_{\cP} &=& \left( \int_{0}^{t \wedge s}
K_{H}(t,u)K_{H}(s,u)du \right)
\langle 1_A, 1_B \rangle_{\cP(\bR^d)} \\
&=&R_H(t,s)\langle 1_A, 1_B \rangle_{\cP(\bR^d)}=\langle 1_{[0,t]
\times A},1_{[0,s] \times B}\rangle_{\cH \cP},
\end{eqnarray*}
i.e. $K_{H}^*$ is an isometry between $(\cE, \langle \cdot,\cdot
\rangle_{\cH \cP})$ and $\cP$. Since $\cH \cP$ is the completion
of $\cE$ with respect to $\langle \cdot, \cdot \rangle_{\cH \cP}$,
this isometry can be extended to $\cH \cP$. We denote this
extension by $K_{\cH \cP}^*$.

\begin{lemma}
\label{isometry-cHcP-L2} $K_{\cH \cP}^* : \cH \cP \rightarrow \cP$
is surjective.
\end{lemma}

\begin{proof}
The proof is similar to the proof of Lemma \ref{isometry-cH-L2},
using the fact that $1_{[0,t] \times A} \in K_{\cH \cP}^{*}(\cH
\cP)$ for all $t \in [0,T], A \in \cB_{b}(\bR^d)$, and $\cE$ is
dense in $\cP$ with respect to $\| \cdot \|_{\cP}$. \qed
\end{proof}

\begin{remark}
{\rm
 Using (\ref{K-H*-as-fract-integr}) and Lemma \ref{isometry-cHcP-L2},
 we can {\em formally} say that
$$\cH \cP=\{\varphi \ \mbox{such that} \ (s,x) \mapsto
s^{-(H-1/2)}I_{T-}^{H-1/2}(u^{H-1/2}\varphi(u,x))(s) \ \mbox{is
in} \ \cP\}.$$ }
\end{remark}

\subsection{The Noise and the Stochastic Integral}

In this subsection, we introduce the noise which is randomly
perturbing the heat equation. This noise is assumed to be
fractional in time and colored in space, with an arbitrary spatial
covariance function $f$. It has been recently considered by other
authors (see e.g. \cite{quersardanyons-tindel06}) in the case when
$f$ is the Riesz kernel and the spatial dimension is $d=1$. The
general definition that we consider in this subsection seems to be
new in the literature.

Let $B=\{B(\varphi); \varphi \in {\cD}((0,T) \times \bR^{d})\}$ be
a zero-mean Gaussian process with covariance
$$E(B(\varphi)B(\psi))=\langle \varphi, \psi \rangle_{\cH \cP}.$$

\noindent Let $H^B$ be the Gaussian space of $B$, i.e. the closed
linear span of $\{B(\varphi); \varphi \in {\cD}((0,T) \times
\bR^{d})\}$ in $L^{2}(\Omega)$. As in subsection
\ref{fractionar-white-noise}, we can define $B_{t}(A)=B(1_{[0,t]
\times A})$ as the $L_{2}(\Omega)$-limit of the Cauchy sequence
$\{B(\varphi)\}_{n}$, where $(\varphi_{n})_{n} \subset \cD((0,T)
\times \bR^d)$ converges pointwise to $1_{[0,t] \times A}$. We
extend this definition by linearity to all elements in $\cE$. A
limiting argument shows that
$$E(B(\varphi)B(\psi))=\langle \varphi, \psi \rangle_{\cH \cP},
\quad \forall \varphi,\psi \in \cE,$$ i.e. $\varphi \mapsto
B(\varphi)$ is an is isometry between $(\cE, \langle \cdot, \cdot
\rangle_{\cH \cP})$ and $H^B$. Since $\cH \cP$ is the completion
of $\cE$ with respect to $\langle \cdot, \cdot \rangle_{\cH}$,
this isometry can be extended to $\cH \cP$, giving us the
stochastic integral with respect to $B$. We will use the notation
$$B(\varphi)=\int_{0}^{T} \int_{\bR^d}\varphi(t,x)B(dt,dx).$$

\begin{remark}
{\rm Similarly to subsection \ref{fractionar-white-noise}, the
transfer operator $K_{\cH \cP}^*$ can be used to explore the
relationship between $B(\varphi)$ and another stochastic integral.
Using Lemma \ref{isometry-cHcP-L2}, we define
\begin{equation}
\label{def-of-M} M(\phi):=B((K_{\cH \cP}^*)^{-1}(\phi)),  \ \ \
\phi \in \cP.
\end{equation} Note that
$$E(M(\phi)M(\eta))=\langle (K_{\cH \cP}^*)^{-1}(\phi),
(K_{\cH \cP}^*)^{-1}(\eta)\rangle_{\cH \cP}=\langle \phi,\eta
\rangle_{\cP},$$ i.e. $M=\{M(\phi); \phi \in \cP\}$ is a Gaussian
noise which is white in time and has spatial covariance function
$f$. This noise has been considered by Dalang in \cite{Da1}. We
use the following notation:
$$M(\phi)=\int_0^T \int_{\bR^d} \phi(t,x)M(dt,dx), \ \ \
\phi \in \cP.$$ Note that $M(\phi)$ is in fact Dalang's stochastic
integral with respect to the noise $M$.

Let $H^M$ be the Gaussian space of $M$, i.e. the closed linear
span of $\{M(\phi); \phi \in \cP\}$ in $L_{2}(\Omega)$. By
(\ref{def-of-M}), it follows that $H^{M}=H^{B}$. The following
diagram summarizes these facts:

\begin{picture}(200,120)(10,10)
\put(20,90){\vector(0,-1){70}} \put(15,95){$\cH \cP$}
\put(110,90){\vector(-1,-1){70}} \put(110,95){$\cP$}
\put(35,100){\vector(1,0){70}} \put(10,10){$H^{B}=H^{M}$}
\put(50,105){$K_{\cH \cP}^{*}$} \put(10,50){$B$} \put(80,50){$M$}
\put(130,70){$B(\varphi)=M(K_{\cH \cP}^*\varphi), \ \ \ \forall
\varphi \in \cH \cP, \ {\rm i.e.}$} \put(130,40){$\int_{0}^{T}
 \int_{\bR^d}\varphi(t,x)B(dt,dx)=\int_{0}^{T}
\int_{\bR^d} (K_{\cH \cP}^*\varphi)
 (t,x)M(dt,dx),\ \ \ \forall \varphi \in \cH \cP.$}
\end{picture}

\noindent In particular, $B(t,A)=\int_{0}^{t}
\int_{A}K_{H}(t,s)M(ds,dy)$. This relationship will not be used in
the present article.

}
\end{remark}

\subsection{The solution of the  Stochastic Heat Equation}

We consider the stochastic heat equation driven by the noise $B$,
written {\em formally} as:

\begin{equation}
\label{heat-eq-HP}u_{t}-\Delta u =  \dot B, \ \ \ \mbox{in} \
(0,T) \times \mathbb{R}^{d}, \ \ \  u(0,\cdot)= 0.
\end{equation}

\noindent As in subsection \ref{fractionar-white-sol}, we let
$G(t,x)$ be given by formula (\ref{fund-sol-heat}), and
$g_{tx}(s,y)=G(t-s,x-y)$. The next theorem is the fundamental
result leading to the necessary and sufficient condition for the
existence of a process-solution and a distribution-solution of
(\ref{heat-eq-HP}).

To state the theorem we need to introduce the following notations:
\begin{eqnarray*}
I_{f,tx}(r,s)&:=&\int_{\bR^d}
\int_{\bR^d}g_{tx}(s,y)f(y-z)g_{tx}(r,z)dy \, dz \\
J_f(u,v,y,z)&:= &\int_{\bR^d} \int_{\bR^d}G(u, y-x) f(x-x')
G(v,z-x') dx \, dx'.
\end{eqnarray*}

\begin{theorem}
\label{HP-psi*Gtilde} Suppose that the spatial covariance function
$f$ satisfies:
\begin{eqnarray}
\label{cond-on-If} A_{f}(2t-s-r)^{-(d-\alpha_f)/2} \leq
I_{f,tx}(r,s) &\leq & B_{f}(2t-r-s)^{-(d-\alpha_f)/2},  \\
\nonumber & & \quad \forall r \in [0,t], \forall s \in
[0,t], \forall t \in [0,T], \forall x \in \bR^d \\
\label{cond-on-Jf} J_f(u,v,y,z) &\leq &
C_{f}(u+v)^{-(d-\alpha_f)/2}, \\
\nonumber & & \forall u \in [0,T], \forall v \in [0,T], \forall y
\in \bR^d, \forall z \in \bR^d
\end{eqnarray}
for some constants $A_{f},B_{f},C_{f}>0$ and $\alpha_f<d$. If
\begin{equation}
\label{cond-HP} H>\frac{d-\alpha_f}{4} \ ,
\end{equation}
then: (a) $g_{tx} \in |\cH \cP|$ for every $(t,x) \in [0,T] \times
\bR^d$; (b) $\eta * \tilde G \in \cH \cP$ for every $\eta \in
\cD((0,T) \times \bR^d)$. Moreover,
\begin{equation}
\label{nec-suf-cond-HP} \|g_{tx}\|_{\cH \cP}<\infty \ \forall
(t,x) \in [0,T] \times \bR^d \quad \mbox{if and only if} \quad
(\ref{cond-HP}) \ \mbox{holds}.
\end{equation}
\end{theorem}

\begin{proof}
(a) We will apply Theorem \ref{criterion-cHcP} to the function
$g_{tx}$. As we noted in the proof of Theorem \ref{H-psi*Gtilde},
$g_{tx}(\cdot,y) \in L_{2}(0,T)$ for every $y \in \bR^d$, and the
map $(\tau,y) \mapsto \cF_{0,T}g_{tx}(\tau,y)$ is measurable. We
calculate
$$\|g_{tx}\|_{\cH \cP}^{2}:=\alpha_H
c_{H}\int_{\bR}\int_{\bR^d}\int_{\bR^d}{\cF}_{0,T}g_{tx}(\tau,x)
f(x-y)\overline{{\cF}_{0,T}g_{tx}(\tau,y)} |\tau|^{-(2H-1)}dy
dxd\tau.$$

\noindent For this, we write
\begin{eqnarray}
\nonumber \|g_{tx}\|_{\cH \cP}^{2} &=& \alpha_H c_H\int_{\bR}
|\tau|^{-(2H-1)} \int_{\bR^d}  \int_{\bR^d}
 \left(\int_{0}^{T}e^{-i\tau s}g_{tx}(s,y)ds \right)f(y-z)
 \left(\int_{0}^{T}e^{i\tau r}g_{tx}(r,y)dr \right) dy \, dz
d\tau \\
\nonumber
 & =& \alpha_H c_H \int_{\bR} |\tau|^{-(2H-1)} \int_{0}^{t}
\int_{0}^{t} e^{-i\tau (s-r)} I_{f,tx}(r,s) dr \, ds \, d\tau
\\
\label{norm-gtx-cHcP}
&=&\alpha_H\int_{0}^{t}\int_{0}^{t}|s-r|^{2H-2}I_{f,tx}(r,s)dr \,
ds,
\end{eqnarray}
where we used Lemma \ref{LemmaA1} (Appendix) for the last equality
and the definition of $I_{f,tx}(r,s)$. Using (\ref{cond-on-If}),
we see that
\begin{eqnarray*}
\|g_{tx}\|_{\cH \cP}^{2} &\leq& \alpha_H B_{f}
\int_{0}^{t}\int_{0}^{t}|s-r|^{2H-2} (2t-r-s)^{-(d-\alpha_f)/2}dr
\, ds \\
\|g_{tx}\|_{\cH \cP}^{2} &\geq& \alpha_H A_{f}
\int_{0}^{t}\int_{0}^{t}|s-r|^{2H-2} (2t-r-s)^{-(d-\alpha_f)/2}dr
\, ds.
\end{eqnarray*}
Relation (\ref{nec-suf-cond-HP}) follows, since the integral above
is finite if and only if $2H>(d-\alpha_f)/2$. In this case, we
have $\|g_{tx}\|_{|\cH \cP|}=\|g_{tx}\|_{\cH \cP}< \infty$ (since
$g_{tx} \geq 0$), and hence $g_{tx} \in |\cH \cP|$.

\vspace{3mm}

(b) We will apply Theorem \ref{criterion-cHcP} to the function
$\varphi=\eta * \tilde G$. By Lemma \ref{eta*G-in-L2},
$\varphi(\cdot,x) \in L_{2}(0,T)$ for every $x \in \bR^d$. We now
calculate
$$\|\varphi\|_{\cH \cP}^{2}:=\alpha_H c_{H} \int_{\bR} \int_{\bR^d}
\int_{\bR^d} {\cF}_{0,T}\varphi(\tau, x)
f(x-y)\overline{{\cF}_{0,T}\varphi(\tau, x')}|\tau|^{-(2H-1)}
dx'dx d\tau.$$

\noindent By (\ref{Fourier-eta*G}), we get
\begin{eqnarray*}
\|\varphi\|^{2}_{\cH \cP} &=&\alpha_H c_H
\int_{\bR}|\tau|^{-(2H-1)} \int _{0}^{T}\int_{0}^{T} e^{-i\tau
(s-r)} \int_{\bR^d} \int_{\bR^d} \int_{0}^{T-s} \int_{0}^{T-r}
\eta (u+s, y) \eta (v+r, z) \\
& & J_f(u,v,y,z)  dv \, du \,dz\, dy  \,dr \, ds \,d\tau \\
& =& \alpha_H \int _{0}^{T}\int_{0}^{T} |s-r|^{2H-2}\int_{\bR^d}
\int_{\bR^d} \int_{0}^{T-s} \int_{0}^{T-r} \eta (u+s, y) \eta
(v+r, z) \\
& & J_f(u,v,y,z)  dv \, du \,dz\, dy \,dr \, ds,
\end{eqnarray*}
where we used Lemma \ref{LemmaA1} (Appendix A) for the last
equality and the definition of $J_{f}(u,v,y,z)$.

Using (\ref{cond-on-Jf}) and the fact that $\eta \in \cD((0,T)
\times \bR^d)$ (and thus is bounded by a constant and its support
is compact), we conclude that
\begin{eqnarray*}
\|\varphi \|_{\cH \cP}^{2} &\leq &\alpha_{H} C_{f} \int_{0}^{T}
\int_{0}^{T} |s-r|^{2H-2} \int_{\bR^d} \int_{\bR^d}
\int_{0}^{T-s} \int_{0}^{T-r} |\eta(u+s,y) \eta(v+r,z)| \\
& & (u+v)^{-(d-\alpha_f)/2}dy \, dz \, dv \, du \, dr \, ds \\
& \leq & \alpha_{H} C_{f} D_{\eta} \int_{0}^{T} \int_{0}^{T}
|s-r|^{2H-2} \int_{0}^{T-s} \int_{0}^{T-r} (u+v)^{-(d-\alpha_f)/2}
dv \, du \, dr \, ds.
\end{eqnarray*}
As in the proof of Theorem \ref{H-psi*Gtilde}.(b), the last
integral is finite if $2H>(d-\alpha_f)/2$.

\qed
\end{proof}

\noindent The next theorem identifies identifies the constant
$\alpha_f$ in the case of some particular covariance functions.

\begin{theorem}
(i) If $f=R_{\alpha}$ with $0<\alpha<d$, then (\ref{cond-on-If})
and (\ref{cond-on-Jf}) hold with $\alpha_f=\alpha$.

\noindent (ii) If $f=B_{\alpha}$ with $\alpha>0$, then
(\ref{cond-on-If}) and (\ref{cond-on-Jf}) hold with $\alpha_f=0$.

\noindent (iii) If $f=G_{\alpha}$ with $\alpha>0$, then
(\ref{cond-on-If}) and (\ref{cond-on-Jf}) hold with $\alpha_f=0$.

\noindent (iv) If $f=P_{\alpha}$ with $\alpha>0$, then
(\ref{cond-on-If}) and (\ref{cond-on-Jf}) hold with $\alpha_f=-1$
\end{theorem}

\begin{remark}
{\rm a) If $f=R_{\alpha}$, condition (\ref{cond-HP}) becomes $H>
\max\{(d-\alpha)/4, 1/2\}$, which does not impose any restrictions
on $d$. For any $d \geq 1$ arbitrary, it suffices to choose
$\alpha$ such that $\max \{d-4H,0\}<\alpha<d$.

\noindent b) If $f=B_{\alpha}$ or $f=G_{\alpha}$, condition
(\ref{cond-HP}) becomes $H> \max\{d/4, 1/2\}$, which forces $d<4$.

\noindent c) If $f=P_{\alpha}$, condition (\ref{cond-HP}) becomes
$H> \max\{(d+1)/4, 1/2\}$, which forces $d<3$. }
\end{remark}

\begin{proof} {\em We begin by examining condition
(\ref{cond-on-If})}. To simplify the notation, we will omit the
index $tx$ in the writing of $I_{f,tx}$. Using the definitions of
$I_{f}$ and $G$, we obtain
\begin{eqnarray}
\nonumber I_f(r,s) &=& \frac{1}{(4 \pi)^d[(t-s)(t-r)]^{d/2}}
\int_{\bR^d} \int_{\bR^d} f(y-z)e^{- \frac{|x-y|^{2}}{4(t-s)}
-\frac{|x-z|^{2}}{4(t-r)}}dy dz \\
\nonumber
 &=&\frac{1}{(2 \pi)^{d}}\int_{\bR^d} \int_{\bR^d}
f(\sqrt{2(t-s)}y'-\sqrt{2(t-r)}z')e^{-\frac{|y'|^2}{2}-
\frac{|z'|^2}{2}}dy' dz' \\
\label{def-of-If} &=& E[f(\sqrt{2(t-s)}Y-\sqrt{2(t-r)}Z)]=E[f(U)].
\end{eqnarray}
Here we used the change of variables $x-y=\sqrt{2(t-s)}y',
x-z=\sqrt{2(t-r)}z'$ and we denoted by $(Y,Z)= (Y_{1}, \ldots
Y_{d}, Z_{1},\ldots , Z_{d}) $ a random vector with independent
$N(0,1)$ components, and
$$U=\sqrt{2(t-s)}Y-\sqrt{2(t-r)}Z.$$

\noindent Note that $U_{i}=\sqrt{2(t-s)}Y_{i}-\sqrt{2(t-r)}Z_{i},
i =1, \ldots, d$ are i.i.d. $N(0,2(2t-s-r))$ random variables.
Then $V_{i}=U_{i}/\sqrt{2(2t-s-r)}, i=1, \ldots, d$ are i.i.d.
$N(0,1)$ random variables and
\begin{equation}
\label{formula-of-U-gtx} |U|^{2}=\sum_{i=1}^{d}U_{i}^{2}=2(2t-s-r)
\sum_{i=1}^{d}V_{i}^2= 2(2t-s-r)W_d,
\end{equation}
 where $W_d$ is a $\chi_d^{2}$
random variable.

\vspace{3mm}

We are now treating separately the four cases:

\vspace{3mm}

(i) In the case of the Riesz kernel,
$f(x)=R_{\alpha}(x)=\gamma_{\alpha,d} |x|^{-(d-\alpha)}$ and
$0<\alpha<d$. Using (\ref{def-of-If}) and
(\ref{formula-of-U-gtx}), the integral $I_{f}(r,s)$ becomes
\begin{eqnarray}
\nonumber
I_{R_{\alpha}}(r,s)&=&\gamma_{\alpha,d}E|U|^{-(d-\alpha)}=\gamma_{\alpha,d}
[2(2t-s-r)]^{-(d-\alpha)/2}E|W_d|^{-(d-\alpha/2)}\\
\label{I-R-alpha} &:=& C_{\alpha,d}(2t-s-r)^{-(d-\alpha)/2},
\end{eqnarray}
 where
$C_{\alpha,d}=
\gamma_{\alpha,d}2^{-(d-\alpha)/2}E|W_d|^{-(d-\alpha)/2}$. This
proves that condition (\ref{cond-on-If}) is satisfied with
$\alpha_f=\alpha$.

\vspace{3mm}

(ii) In the case of the Bessel kernel,
$f(x)=B_{\alpha}(x)=\gamma_{\alpha}'\int_{0}^{\infty }
w^{(\alpha-d)/2-1} e^{-w}e^{-|x|^{2}/4w} dw$ and $\alpha>0$. Using
(\ref{def-of-If}) and (\ref{formula-of-U-gtx}), the integral
$I_{f}(r,s)$ becomes
\begin{eqnarray*}
I_{B_{\alpha}}(r,s)&=&\gamma'_{\alpha}
\int_{0}^{\infty}w^{-(\alpha-d)/2-1}
e^{-w} E[e^{-|U|^2/(4w)}]dw \\
&=& \gamma'_{\alpha}\int_{0}^{\infty}w^{(\alpha-d)/2-1} e^{-w}
E\left[\exp\left(-\frac{2t-r-s}{2w} W_d\right)\right]dw \\
 &=& \gamma'_{\alpha}\int_{0}^{\infty}w^{(\alpha-d)/2-1} e^{-w}
\left(1+\frac{2t-r-s}{w} \right)^{-d/2}dw
\end{eqnarray*}
where we used the fact that $E[e^{-cW_d}]=(1+2c)^{-d/2}$ for any
$c>0$. Note that
$$ \left(\frac{2t-r-s}{w}\right)^{d/2} \leq
\left(1+\frac{2t-r-s}{w}\right)^{d/2} \leq
C_{d}\left[1+\left(\frac{2t-r-s}{w}\right)^{d/2} \right]$$ where
$C_{d}=2^{d/2-1}$. Hence
\begin{equation}
\label{double-ineq} \frac{1}{2C_{d}} (2t-r-s)^{-d/2} \leq
\frac{1}{C_{d}} \frac{w^{d/2}}{w^{d/2}+(2t-s-r)^{d/2}}\leq
\left(1+\frac{2t-r-s}{w}\right)^{-d/2} \leq w^{d/2}(2t-s-r)^{-d/2}
\end{equation}
where for the first inequality we used the fact that
$a/(a+x)>1/(2x)$ if $x$ is small enough and $a>0$.

We conclude that
\begin{eqnarray*}
I_{B_{\alpha}}(r,s) &\leq &\gamma'_{\alpha}\int_{0}^{\infty}
w^{(\alpha-d)/2-1} e^{-w} w^{d/2}(2t-r-s)^{-d/2}dw \leq
\gamma_{\alpha}'\Gamma\left(\frac{\alpha}{2}\right)
(2t-r-s)^{-d/2} \\
I_{B_{\alpha}}(r,s) & \geq &
\frac{\gamma_{\alpha}'}{2C_d}\int_{0}^{1} w^{(\alpha-d)/2-1}
e^{-w} (2t-r-s)^{-d/2}dw \geq \frac{\gamma_{\alpha}'}{2C_d}
\left(\int_{0}^{1} w^{\alpha/2-1} e^{-w} dw \right)(2t-r-s)^{-d/2}
,
\end{eqnarray*}
i.e. condition (\ref{cond-on-If}) is satisfied with $\alpha_f=0$.

\vspace{3mm}

(iii) In the case of the heat kernel,
$f(x)=G_{\alpha}(x)=\gamma_{\alpha,d}''e^{-|x|^{2}/(4\alpha)}$ and
$\alpha>0$. Using (\ref{def-of-If}) and (\ref{formula-of-U-gtx}),
the integral $I_{f}(r,s)$ becomes: $$I_{G_{\alpha}} (r,s)= \gamma
_{\alpha,d}'' E[e^{-|U|^{2}/(4\alpha)}] =\gamma _{\alpha,d}'' E
\left[\exp \left( -\frac{2t-s-r}{2 \alpha}W_{d} \right)\right] =
\gamma _{\alpha,d}'' \left( 1+ \frac{2t-r-s}{\alpha
}\right)^{-d/2}.$$

\noindent Using (\ref{double-ineq}), we obtain that
$$\frac{ \gamma_{\alpha,d}''}{2C_{d}}(2t-r-s)^{-d/2} \leq
I_{G_{\alpha}} (r,s) \leq
 \gamma_{\alpha,d}''\alpha^{d/2} (2t-r-s)^{-d/2},$$
i.e. condition (\ref{cond-on-If}) is satisfied with $\alpha_f=0$.

\vspace{3mm}

(iv) In the case of the Poisson kernel,
$f(x)=P_{\alpha}(x)=\gamma_{\alpha,d}'''
(|x|^{2}+\alpha^{2})^{-(d+1)/2}$ and $\alpha>0$. Using
(\ref{def-of-If}) and (\ref{formula-of-U-gtx}), the integral
$I_{f}(r,s)$ becomes: $$I_{P_{\alpha}}(r,s)= \gamma _{\alpha,d}'''
E\left[(|U|^2+\alpha^2)^{-(d+1)/2}\right]=\gamma _{\alpha,d}'''
E\left|2(2t-r-s)W_d+\alpha^2\right|^{-(d+1)/2}.$$ Using the fact
that
$$A_{d}[(2t-r-s)W_d]^{-(d+1)/2} \leq [2(2t-r-s)W_d+\alpha^2]^{-(d+1)/2}
\leq  B_{d} [(2t-r-s)W_d]^{-(d+1)/2}$$ for some constants
$A_d,B_d>0$, we conclude that
$$A_{d}E|W_d|^{-(d+1)/2}(2t-r-s)^{-(d+1)/2} \leq
I_{P_{\alpha}} (r,s) \leq B_d
E|W_d|^{-(d+1)/2}(2t-r-s)^{-(d+1)/2},$$ i.e. condition
(\ref{cond-on-If}) is satisfied with $\alpha_f=-1$.

\vspace{3mm}

{\em We continue by examining condition (\ref{cond-on-Jf})}. Using
the definitions of $J_f$ and $G$, we obtain that
\begin{eqnarray}
\nonumber J_f(u,v,y,z)&=&\frac{1}{(4
\pi)^{d}(uv)^{d/2}}\int_{\bR^d} \int_{\bR^d}
f(x-x')e^{-\frac{|x-y|^{2}}{4u}
-\frac{|x'-z|^{2}}{4v}}dx dx' \\
\nonumber &=&
\frac{1}{(2\pi)^{d}}\int_{\bR^d}\int_{\bR^d}f(y-z+\sqrt{2}(\sqrt{u}a+
\sqrt{v}a'))e^{-\frac{|a|^{2}}{2}-\frac{|a'|^{2}}{2}}da da' \\
\label{J-f-calcul} &=& E[f(y-z+\sqrt{2}
(\sqrt{u}Y+\sqrt{v}Z))]=E[f(y-z+U)].
\end{eqnarray}

\noindent Here we used the change of variables $x-y=\sqrt{2u}a$
and $x'-z=\sqrt{2v}a'$, we denoted by $(Y,Z)= (Y_{1}, \ldots
Y_{d}, Z_{1},\ldots , Z_{d}) $ a random vector with independent
$N(0,1)$ components, and
$$U=\sqrt{2u}Y-\sqrt{2v}Z.$$

\noindent Note that $U_{i}=\sqrt{2u}Y_i-\sqrt{2v}Z_i,
i=1,\ldots,d$ are i.i.d. $N(0,2(u+v))$ random variables. Then
$V_{i}=U_{i}/\sqrt{2(u+v)},i=1, \ldots,d$ are i.i.d $N(0,1)$
random variables and
\begin{equation}
\label{norm-y-z+U}
|y-z+U|^2=\sum_{i=1}^{d}(y_i-z_i+U_i)^2=2(u+v)\sum_{i=1}^{d}
(\mu_i+V_{i})^{2}=2(u+v)\sum_{i=1}^{d}T_{i}^{2},
\end{equation}
where $\mu_{i}=(y_i-z_i)/\sqrt{2(u+v)}$ and $T_{i}=\mu_i+V_i$ is
$N(\mu_i,1)$-distributed. It is known that (see e.g. p. 132,
\cite{johnson-kotz70})
\begin{equation}
\label{non-central-chi2}
\sum_{i=1}^{d}T_i^{2}\stackrel{d}{=}W_{d-1}+S^2,
\end{equation}
where $W_{d-1}$ and $S$ are independent random variables with
distributions $\chi_{d-1}^2$, respectively
$N(\sqrt{\sum_{i=1}^{d}\mu_i^2},1)$.

We are now treating separately the four cases:

\vspace{3mm}

(i) In the case of the Riesz kernel,
$f(x)=\gamma_{\alpha,d}|x|^{-(d-\alpha)}$. Using
(\ref{J-f-calcul}), (\ref{norm-y-z+U}), and
(\ref{non-central-chi2}), the integral $J_{f}(u,v,y,z)$ becomes:
\begin{eqnarray*}
J_{R_{\alpha}}(u,v,y,z)&=&\gamma_{\alpha,d}E|y-z+U|^{-(d-\alpha)}=
\gamma_{\alpha,d} [2(u+v)]^{-(d-\alpha)/2}
E\left|\sum_{i=1}^{d}T_i^2\right|^{-(d-\alpha)/2} \\
 &=& \gamma_{\alpha,d} [2(u+v)]^{-(d-\alpha)/2}
E\left|W_{d-1}+S^2\right|^{-(d-\alpha)/2} \\
&\leq & D_{\alpha,d} (u+v) ^{-(d-\alpha )/2},
\end{eqnarray*}
where $D_{\alpha,d}=
\gamma_{\alpha,d}2^{-(d-\alpha)/2}E|W_{d-1}|^{-(d-\alpha)/2}$,
i.e. condition (\ref{cond-on-Jf}) is satisfied with
$\alpha_f=\alpha$.

\vspace{3mm}

(ii) In the case of the Bessel kernel,
$f(x)=B_{\alpha}(x)=\gamma_{\alpha}'\int_{0}^{\infty }
w^{(\alpha-d)/2-1} e^{-w}e^{-|x|^{2}/4w} dw$ and $\alpha>0$. Using
(\ref{J-f-calcul}), (\ref{norm-y-z+U}), and
(\ref{non-central-chi2}), the integral $I_{f}(r,s)$ becomes
\begin{eqnarray*}
J_{B_{\alpha}}(u,v,y,z)&=&\gamma_{\alpha}'\int_{0}^{\infty}w^{(\alpha-d)/2-1}
e^{-w}E[e^{-|y-z+U|^{2}/(4w)}]dw \\
&=&\gamma_{\alpha}'\int_{0}^{\infty}w^{(\alpha-d)/2-1}
e^{-w}E\left[\exp \left(-\frac{u+v}{2w}\sum_{i=1}^{d}T_{i}^2
\right) \right]dw \\
&=&\gamma_{\alpha}'\int_{0}^{\infty}w^{(\alpha-d)/2-1}
e^{-w}E\left\{\exp \left[-\frac{u+v}{2w}(W_{d-1}+S^2) \right]
\right\}dw
\end{eqnarray*}

\noindent Note that
\begin{equation}
\label{mgf-W-S} E[e^{-c(W_{d-1}+S^2)}] \leq \left(
1+2c\right)^{-d/2}, \quad \forall c>0.
\end{equation}

\noindent This follows by the independence of $W_{d-1}$ and $S$,
the fact that $E(e^{-cW_{d-1}})=(1+2c)^{-(d-1)/2}$, and
$$E\left( e^{-cS^{2} } \right) = \frac{1}{\sqrt{1+2c}} \exp
\left \{-\frac{\sum_{i=1}^{d}\mu_{i}^{2}}{2} \cdot
\frac{2+2c}{1+2c} \right\} \leq (1+2c)^{-1/2}, \quad \forall c>0$$
(recall that $S$ has $N(\sqrt{\sum_{i=1}^{d}\mu_{i}^{2}},1)$
distribution). Therefore
\begin{eqnarray*}
J_{B_{\alpha}}(u,v,y,z)&\leq & \gamma _{\alpha }'\int_{0}^{
\infty}w^{(\alpha-d)/2-1} e^{-w}\left(
1+\frac{u+v}{w}\right)^{-d/2}dvdu \\
& = &\gamma _{\alpha}' \int_{0}^{\infty } w^{\alpha /2-1}
e^{-w} (w+u+v) ^{-d/2} dw \\
& \leq & \gamma_{\alpha}'\Gamma(\alpha/2) (u+v)^{-d/2},
\end{eqnarray*}
i.e. condition (\ref{cond-on-Jf}) is satisfied with $\alpha_f=0$.

\vspace{3mm}

(iii) In the case of the heat kernel,
$f(x)=G_{\alpha}(x)=\gamma_{\alpha,d}''e^{-|x|^{2}/(4\alpha)}$ and
$\alpha>0$. Using (\ref{J-f-calcul}), (\ref{norm-y-z+U}), and
(\ref{non-central-chi2}), the integral $J_{f}(u,v,y,z)$ becomes
\begin{eqnarray*}
J_{G_{\alpha}}(u,v,y,z)&=& \gamma_{\alpha,d}''
E[e^{-|U|^{2}/(4\alpha)}]=\gamma_{\alpha,d}'' E\left\{\exp
\left[-\frac{u+v}{2\alpha}(W_{d-1}+S^{2}) \right]
\right\} \\
& \leq & \gamma _{\alpha,d}'' \left( 1+ \frac{u+v}{\alpha }
\right) ^{-d/2} \leq \gamma_{\alpha,d}''\alpha^{d/2} (u+v) ^{-d/2}
\end{eqnarray*}
where we used (\ref{mgf-W-S}) for the first inequality. This
proves that condition (\ref{cond-on-Jf}) is satisfied with
$\alpha_f=0$.

\vspace{3mm}

(iv) In the case of the Poisson kernel,
$f(x)=P_{\alpha}(x)=\gamma_{\alpha,d}'''
(|x|^{2}+\alpha^{2})^{-(d+1)/2}$ and $\alpha>0$. Using
(\ref{J-f-calcul}), (\ref{norm-y-z+U}), and
(\ref{non-central-chi2}), the integral $J_{f}(u,v,y,z)$ becomes
\begin{eqnarray*}
J_{P_{\alpha}}(u,v,y,z)&=&
\gamma_{\alpha,d}'''E\left[(|y-z+U|^2+\alpha^2)^{-(d+1)/2}\right]
= \gamma_{\alpha,d}'''E\left|2(u+v)\sum_{i=1}^{d}T_{i}^2+
\alpha^2\right|^{-(d+1)/2} \\
&=&\gamma_{\alpha,d}'''E\left|2(u+v)(W_{d-1}+S^2)+
\alpha^2\right|^{-(d+1)/2} \leq \gamma_{\alpha,d}'''
[2(u+v)]^{-(d+1)/2}E|W_{d-1}|^{-(d+1)/2},
\end{eqnarray*}
i.e. condition (\ref{cond-on-Jf}) is satisfied with $\alpha_f=-1$.
\qed
\end{proof}

\vspace{3mm}

Under the conditions of Theorem \ref{HP-psi*Gtilde}), $B(g_{tx})$
and $B(\eta * \tilde G)$ are well-defined for every $(t,x)$,
respectively for every $\eta \in {\cD}((0,T) \times \bR^{d})$, and
we can introduce the following definition:

\begin{definition}
{\rm a) The process $\{u(t,x);t \in [0,T], x \in \bR^{d}\}$
defined by
\begin{equation}
\label{solution-process-HP} u(t,x):=B(g_{tx})=\int_{0}^{T}
\int_{\bR^{d}} G(t-s,x-y)B(ds,dy)
\end{equation}
is called the {\bf process solution} of the stochastic heat
equation (\ref{heat-eq-HP}), with vanishing initial conditions.

b) The process $\{u(\eta); \eta \in {\cD}((0,T) \times
\mathbb{R}^{d})\}$ defined by
$$u(\eta):  = B(\eta * \tilde G)=\int_{0}^{T} \int_{\bR^{d}}
\left( \int_{0}^{\infty} \int_{\bR^{d}}\eta(t+s,x+y)G(s,y)dyds
\right) B(dt,dx)$$ is called the {\bf distribution-valued
solution} of the stochastic heat equation (\ref{heat-eq-HP}), with
vanishing initial conditions.}
\end{definition}

\begin{lemma}
\label{L2-cont-u} The process $\{u(t,x);t \in [0,T], x \in
\bR^{d}\}$ is $L^{2}(\Omega)$-continuous.
\end{lemma}

\begin{proof}
The proof is identical to the proof of Lemma \ref{L2-cont-v},
based on the continuity of the function $G(t,x)$ with respect to
each of its arguments, and the fact that $\|g_{tx}\|_{\cH
\cP}<\infty$, which is a consequence of Theorem
\ref{HP-psi*Gtilde}.(a). \qed
\end{proof}

\vspace{3mm}

The next theorem can be viewed as a counterpart of the result
obtained by Maslowski and Nualart (see Example 3.5, \cite{MaNu} in
the case $m=1$, $L_1=\Delta$, $f=0$, $\Phi=1$).

\begin{theorem}
\label{existence-theorem-HP} Suppose that the spatial covariance
function $f$ satisfies (\ref{cond-on-If}) and (\ref{cond-on-Jf}).
Let $\{u(\eta); \eta \in {\cD}((0,T) \times \mathbb{R}^{d})\}$ be
the distribution-valued solution of the stochastic heat equation
(\ref{heat-eq-HP}).

In order that there exists a jointly measurable and locally
mean-square bounded process $X=\{X(t,x);t \in [0,T],x \in
\mathbb{R}^{d}\}$ such that
$$u(\eta)=\int_{0}^{T} \int_{\bR^{d}} X(t,x)\eta(t,x)
dxdt \ \ \ \forall \eta \in {\cD}((0,T) \times \bR^{d}) \ \ \
a.s.$$ it is necessary and sufficient that (\ref{cond-HP}) holds.
In this case, $X$ is a modification of the process $u=\{u(t,x);t
\in [0,T],x \in \bR^{d}\}$ defined by (\ref{solution-process-HP}).
\end{theorem}

\begin{proof}
The proof is omitted since it is identical to the proof of Theorem
\ref{existence-theorem-H}, using $\|\cdot \|_{\cH \cP}$ instead of
$\|\cdot \|_{\cH}$, and relation (\ref{nec-suf-cond-HP}) instead
of (\ref{nec-suf-cond-H}). \qed
\end{proof}

\appendix

\section{An Auxiliary Lemma}
The following result is the analogue of Lemma 1, p.116,
\cite{stein70}, for functions on bounded domains. It plays a
crucial role in the present paper. For our purposes, it is stated
only for $d=1$, but it can be easily generalized to $d \geq 2$.

\begin{lemma}
\label{LemmaA1} Let $0<\alpha<1$ be arbitrary. (a) For every
$\varphi \in L_2(a,b)$, we have
$$\int_{a}^{b}|t|^{-(1-\alpha)} \varphi(t)dt=q_{\alpha}\int_{\bR}|\tau|^{-\alpha}
{\cF}_{a,b}\varphi(\tau)d\tau$$ where
$q_{\alpha}=(2^{1-\alpha}\pi^{1/2})^{-1}\Gamma(\alpha/2)/
\Gamma((1-\alpha)/2)$. (b) For every $\varphi,\psi \in L_2(a,b)$,
we have
$$\int_{a}^{b}
\int_{a}^{b}\varphi(u)|u-v|^{-(1-\alpha)}\psi(v)dvdu=q_{\alpha}\int_{\bR}
|\tau|^{-\alpha} {\cF}_{a,b} \varphi (\tau) \overline{{\cF}_{a,b}
\psi (\tau)} d\tau$$
\end{lemma}

\begin{remark}
{\rm Note that $q_{\alpha}=1/\gamma_{\alpha,1}$ where
$\gamma_{\alpha,d}$ is the constant defined in Example
\ref{riesz}. }
\end{remark}

\begin{proof} (a) We use the fact that
\begin{equation}
\label{Fourier-expon} \int_{\bR}e^{-i \tau t}e^{-\pi \delta
|\tau|^{2}} d \tau=\delta^{-1/2}e^{-|t|^{2}/(4 \pi\delta)}, \ \ \
\forall \delta>0.
\end{equation}

\noindent Using the definition of ${\cF}_{a,b} \varphi$, Fubini's
theorem and (\ref{Fourier-expon}), we have
$$\int_{\bR}e^{-\pi
\delta |\tau|^2} {\cF}_{a,b} \varphi(\tau)d \tau=
\int_{a}^{b}\left(\int_{\bR}e^{-i\tau t}e^{-\pi \delta |\tau|^2}d
\tau \right)\varphi(t) dt = \delta^{-1/2}\int_{a}^{b}
e^{-|t|^{2}/(4\pi\delta)} \varphi(t)dt.$$

\noindent Multiply by $\delta^{\alpha/2-1}$ and integrate with
respect to $\delta>0$. Using Fubini's theorem, we get
$$\int_{\bR} \left(\int_{0}^{\infty} \delta^{\alpha/2-1} e^{-\pi \delta |\tau|^2}
d\delta \right){\cF}_{a,b} \varphi(\tau)d \tau=\int_{a}^{b}
\left(\int_{0}^{\infty} \delta^{-(1-\alpha)/2-1}e^{-\pi
|t|^{2}/(4\pi\delta)} d \delta \right) \varphi(t)dt.$$ Using the
change of variable $1/\delta=u$ for the inner integral on the
right hand side, and the definition of the Gamma function for
evaluating both inner integrals, we get the conclusion.

(b) Note that for every $u \in [a,b]$, $$\int_{a}^{b}
|u-v|^{-(1-\alpha)} \psi(v)dv= \int_{u-b}^{u-a}|w|^{-(1-\alpha)}
\psi(u-w)dw =
q_{\alpha}\int_{\bR}|\tau|^{-\alpha}{\cF}_{u-b,u-a}(\psi_{u})^{\verb2~2}(\tau)d
\tau,$$
 where we used the result in (a) for the last
equality. Now,
$${\cF}_{u-b,u-a}(\psi_{u})^{\verb2~2}(\tau)=\int_{u-b}^{u-a}e^{-\tau t}
\psi(u-v)dv=\int_{a}^{b}e^{-i\tau(u-w)}\psi(w)dw=e^{-i\tau u}
\overline{{\cF}_{a,b} \psi(\tau)}.$$

\noindent Using Fubini's theorem
\begin{eqnarray*}
\lefteqn{\int_{a}^{b}\varphi(u)\int_{a}^{b}|u-v|^{-(1-\alpha)}\psi(v)
dv du=q_{\alpha}\int_{a}^{b}
\varphi(u)\int_{\bR}|\tau|^{-\alpha}e^{-i\tau u}
\overline{{\cF}_{a,b} \psi(\tau)}d \tau du =} \\
&& q_{\alpha}\int_{\bR}\left( \int_{a}^{b}\varphi(u)e^{-i\tau u}du
\right) |\tau|^{-\alpha} \overline{{\cF}_{a,b} \psi(\tau)} d\tau=
q_{\alpha}\int_{\bR}
 |\tau|^{-\alpha} {\cF}_{a,b} \varphi(\tau)  \overline{{\cF}_{a,b}
\psi(\tau)}  d\tau.
\end{eqnarray*}
\qed \end{proof}

\section{Proof of (\ref{Theorem11-Dalang-claim})}

\noindent Let $\eta_n(t,x)=\lambda_n(t-t_0)
\psi_n(x-x_0):=\alpha_n(t)\beta_n(x)$. Then
\begin{eqnarray*}
\cF_{0,T}(\eta_n* \tilde G)(\tau,x)&=& 
 \int_{0}^{T}e^{-i \tau t} \left(\int_{-\infty}^{0}
\int_{\bR^d} \alpha_n(t-s) \beta_n(x-y)\tilde G(s,y)dy \, ds
\right)
dt \\
&=& \int_{\bR^d}\beta_n(x-y)  \int_{-\infty}^{0}e^{-i \tau s}
\tilde G(s,y) \left( \int_{0}^{T} e^{-i \tau (t-s)} \alpha_n(t-s)
dt \right) ds \, dy.
\end{eqnarray*}

\noindent Since ${\rm supp} \ \alpha_n \subset (t_0,t_0+T/n)$, we
obtain that
$$\int_{0}^{T} e^{-i \tau (t-s)} \alpha_n(t-s) dt
=\int_{-s}^{T-s}e^{-i \tau u}\alpha_n(u) du= \left\{
\begin{array}{ll}
\cF_{t_0,t_0+T/n}\alpha_n(\tau) & \mbox{if $-s<t_0$} \\
\cF_{-s,t_0+T/n}\alpha_n(\tau) & \mbox{if $t_0<-s<t_0+T/n$} \\
0 & \mbox{if $-s>t_0+T/n$}
\end{array}
\right.$$ and hence
$$\cF_{0,T}(\eta_n* \tilde G)(\tau,x)=\cF_{t_0,t_0+T/n}\alpha_n(\tau)
\int_{\bR^d}\beta_n(x-y)  \int_{-t_0}^{0}e^{-i \tau s} \tilde
G(s,y) ds dy+ $$ $$\int_{\bR^d}\beta_n(x-y)
\int_{-t_0-T/n}^{-t_0}e^{-i \tau s} \tilde G(s,y)
\cF_{-s,t_0+T/n}\alpha_n(\tau) ds
dy:=A_{n}(\tau,x)+B_{n}(\tau,x).$$

\noindent Note that $\lim_n B_{n}(\tau,x)=0$. Whereas for
$A_{n}(\tau,x)$, we have
$$\cF_{t_0,t_0+T/n}\alpha_n(\tau)=\int_{t_0}^{t_0+T/n}e^{-i \tau t}
\lambda_n(t-t_0)dt=e^{-i \tau t_0}\int_{0}^{T/n}e^{-i \tau
u}\lambda_n(u)du$$ $$=e^{-i \tau t_0} \cF_{0,T} \lambda(\tau/n)
\to e^{-i \tau t_0}, \quad {\rm as} \ n \to \infty$$ and
\begin{eqnarray*}
\int_{\bR^d}\beta_n(x-y)  \int_{-t_0}^{0}e^{-i \tau s} \tilde
G(s,y) ds dy&=& \int_{\bR^d}\psi_n(x-y-x_0)\cF_{-t_0,0} \tilde
G(\tau,y)dy\\
&=& \int_{\bR^d}\psi_n(z-x_0)\cF_{-t_0,0} \tilde G(\tau,x-z)dz \\
&\to &\cF_{-t_0,0} \tilde G(\tau,x-x_0), \quad {\rm as} \ n \to
\infty
\end{eqnarray*}
 by Lebesgue differentiation theorem (see Exercise 2, Chapter 7,
\cite{wheeden-zygmund77}). Therefore
$$\lim_n A_{n}(\tau,x)=e^{-i \tau t_0}\cF_{-t_0,0} \tilde
G(\tau,x-x_0)= \cF_{0,T}g_{t_0 x_0}(\tau,x).$$

\section{A Result about the Space $\cP(\bR^d)$}

The next result is embedded in Theorem 3, \cite{Da1}. We have used
this result in the proof of Lemma \ref{indic-approx-Fourier2}. We
include its proof for the sake of completeness.

\begin{lemma}
\label{LemmaC1} If $\varphi \in  \cS(\bR^d)$ and
$\int_{\bR^d}|\cF_{2}\varphi(\xi)|^2 \mu(d\xi)<\infty$, then
$\varphi \in \cP(\bR^d)$.
\end{lemma}

\begin{proof}
Let $\eta \in \cD(\bR^d)$ be such that $\eta>0$ and
$\int_{\bR^d}\eta(x)dx=1$. Define $\eta_n(x)=n^d \eta(nx)$ and
$\varphi_{n}=\varphi*\eta_n \in \cS(\bR^d)$. We have $\varphi_n
\in |\cP(\bR^d)| \subset \cP(\bR^d)$, since
$$\int_{\bR^d}\int_{\bR^d} |\varphi_n(x)|f(x-y)|\varphi_{n}(y)|dy
dx=\int_{\bR^d}f(z)(|\varphi_n|*|\tilde \varphi_n|)(z)dz<\infty,$$
by Leibnitz's formula (see p. 13, \cite{Da1}). Note that
$$\|\varphi_n-\varphi \|_{\cP(\bR^d)}^{2}=\int_{\bR^d}|\cF_{2}\varphi_n(\xi)-
\cF_{2}\varphi(\xi)|^2
\mu(d\xi)=\int_{\bR^d}|\cF_{2}\eta_{n}(\xi)-1|^{2}|\cF_{2}\varphi(\xi)|^2
\mu(d\xi) \to 0,$$ where we used the dominated convergence
theorem, and the fact that $\cF_{2}\eta(\xi)=\cF_{2}\eta(\xi/n)
\to 1$ and $|\cF_{2}\eta(\xi)| \leq 1$ for all $n$. The conclusion
follows.

\qed
\end{proof}

\end{document}